\documentclass[a4paper,reqno,twoside]{amsart}
\usepackage{bbm}
\usepackage[left=3.5cm,right=3.5cm,top=3.5cm,bottom=3.5cm,footskip=1.5cm,headsep=1cm,bindingoffset=0.cm,]{geometry}

\usepackage[utf8]{inputenc}
\usepackage[english]{babel}
\usepackage[T1]{fontenc} 
\usepackage{lmodern} 
\rmfamily 
\DeclareFontShape{T1}{lmr}{b}{sc}{<->ssub*cmr/bx/sc}{}
\DeclareFontShape{T1}{lmr}{bx}{sc}{<->ssub*cmr/bx/sc}{}

\numberwithin{equation}{section}
\usepackage{amsmath}
\usepackage{amssymb}
\usepackage{amsthm}
\usepackage{amsfonts}
\usepackage{mathtools}

\usepackage{mathdots}
\usepackage{caption}
\usepackage{subcaption}

\usepackage{dsfont} 
\usepackage[format=hang]{caption}
\usepackage[normalem]{ulem}

\usepackage{standalone}
\usepackage{graphicx}
\usepackage{imakeidx}
\makeindex
\usepackage{setspace}
\usepackage{appendix}
\usepackage{pdfpages}
\usepackage{upgreek}
\usepackage{tabulary}
\usepackage{booktabs}
\usepackage{extarrows}
\usepackage{tabularx}
\usepackage{float}
\restylefloat{table}
\usepackage{enumitem}
\usepackage{csquotes}
\usepackage{bm}
\usepackage{orcidlink}
\usepackage{lipsum}

\usepackage{tikz}
\usepackage{tikz-layers}
\usepackage{tikz-cd}
\usepackage[arrow, matrix, curve]{xy}
\usetikzlibrary{matrix,positioning,decorations.pathreplacing,calc}
\usetikzlibrary{
    decorations.text,%
    decorations.markings,%
    shadows}
\usetikzlibrary{calc,intersections}
\usepackage{adjustbox}
\usepackage{graphicx}

\usepackage[
    style=numeric,
    sorting=nty,
    maxnames=99,
    maxalphanames=5,
    natbib=true,
    backend=biber,
    sortcites]{biblatex}

\DeclareNameAlias{default}{family-given}

\AtEveryBibitem{
    \clearfield{url}
    \clearfield{issn}
    \clearfield{isbn}
    \clearfield{urldate}

    \ifentrytype{book}{
        \clearfield{pages}}{%
    }
}

\usepackage{xargs}                      
\usepackage[prependcaption,textsize=tiny,textwidth=3cm]{todonotes}
\newcommandx{\unsure}[2][1=]{\todo[linecolor=red,backgroundcolor=red!25,bordercolor=red,#1]{#2}}
\newcommandx{\change}[2][1=]{\todo[linecolor=blue,backgroundcolor=blue!25,bordercolor=blue,#1]{#2}}
\newcommandx{\info}[2][1=]{\todo[linecolor=OliveGreen,backgroundcolor=OliveGreen!25,bordercolor=OliveGreen,#1]{#2}}
\newcommandx{\improvement}[2][1=]{\todo[linecolor=black,backgroundcolor=black!25,bordercolor=black,#1]{#2}}
\newcommandx{\thiswillnotshow}[2][1=]{\todo[disable,#1]{#2}}

\usepackage[noabbrev,capitalize]{cleveref}
\crefname{proposition}{Proposition}{Propositions}
\crefname{equation}{}{}

\newtheorem{theorem}{Theorem}[section]
\newtheorem{lemma}[theorem]{Lemma}
\newtheorem{proposition}[theorem]{Proposition}
\newtheorem{corollary}[theorem]{Corollary}

\theoremstyle{definition}
\newtheorem{definition}[theorem]{Definition}

\newtheorem{remark}[theorem]{Remark}

\crefname{assumption}{Assumption}{Assumptions}
\crefname{definition}{Definition}{Definitions}
\crefname{corollary}{Corollary}{Corollaries}
\crefname{enumi}{item}{items}

\creflabelformat{subfigure}{#2\textsc{#1}#3}

\usepackage[final]{microtype}

\makeatletter
\newsavebox\myboxA
\newsavebox\myboxB
\newlength\mylenA

\newcommand*\xoverline[2][0.75]{%
  \sbox{\myboxA}{$\m@th#2$}%
  \setbox\myboxB\null
  \ht\myboxB=\ht\myboxA%
  \dp\myboxB=\dp\myboxA%
  \wd\myboxB=#1\wd\myboxA
  \sbox\myboxB{$\m@th\overline{\copy\myboxB}$}
  \setlength\mylenA{\the\wd\myboxA}
  \addtolength\mylenA{-\the\wd\myboxB}%
  \ifdim\wd\myboxB<\wd\myboxA%
    \rlap{\hskip 0.5\mylenA\usebox\myboxB}{\usebox\myboxA}%
  \else
    \hskip -0.5\mylenA\rlap{\usebox\myboxA}{\hskip 0.5\mylenA\usebox\myboxB}%
  \fi}
\makeatother

\usepackage{fancyhdr}

\fancypagestyle{plain}{%
    \fancyhead[C]{}
    \fancyfoot[C]{}
    
    }
\fancypagestyle{nosection}{%
    \fancyhead[CE]{}
    \fancyhead[CO]{}
    \fancyfoot[CE,CO]{\thepage}
}

\usepackage{tikz}

\usetikzlibrary{matrix} 
\usetikzlibrary{arrows,arrows.meta,bending} 

\usepackage{nicematrix}

\usetikzlibrary{hobby}

\usepackage{calc}

\DeclareMathOperator{\Z}{\mathbb{Z}}

\DeclareMathOperator{\R}{\mathbb{R}}
\DeclareMathOperator{\C}{\mathbb{C}}

\renewcommand{\tilde}{\widetilde}

\renewcommand{\qedsymbol}{$\blacksquare$}

\renewcommand{\epsilon}{\varepsilon}

\let\emptyset\varnothing

\renewcommand{\tilde}{\widetilde}

\usepackage{tcolorbox}
\usetikzlibrary{tikzmark,calc}

\title[Off-diagonal decay of capacitance matrices via Brownian motion]{Off-diagonal decay of capacitance matrices via Brownian motion}

\begin{document}

\author[H. Ammari]{Habib Ammari \,\orcidlink{0000-0001-7278-4877}}
\address{\parbox{\linewidth}{Habib Ammari\\
 ETH Z\"urich, Department of Mathematics, Rämistrasse 101, 8092 Z\"urich, Switzerland \href{http://orcid.org/0000-0001-7278-4877}{orcid.org/0000-0001-7278-4877}}.}
 \email{habib.ammari@math.ethz.ch}

\author[C. Thalhammer]{Clemens Thalhammer \,\orcidlink{0009-0006-7218-260X}}
\address{\parbox{\linewidth}{Clemens Thalhammer\\
 ETH Z\"urich, Department of Mathematics, 
 R\"amistrasse 101, 8092 Z\"urich, Switzerland, 
 \href{https://orcid.org/0009-0006-7218-260X}{orcid.org/0009-0006-7218-260X}}.}
\email{clemens.thalhammer@sam.math.ethz.ch}
\thanks{}

 \begin{abstract}
    We establish upper and lower bounds for the off-diagonal entries of the capacitance matrix of infinite systems of spherical resonators in chain, screen and crystal geometries. Let $d$ denote the distance between two resonator centres. For chains contained in a cylinder of bounded radius and with bounded axial gaps, the decay is of order $d^{-1}\log^{-2}(2+d)$. For screens contained in a slab and relatively dense in the tangential directions, it is of order $d^{-3}$. For relatively dense three-dimensional configurations, we obtain exponential upper and lower bounds. No periodicity is assumed. Our proofs are based on the classical connection between the Dirichlet problem and Brownian motion, which enables us to relate the decay estimates to hitting probabilities of Brownian motion.
 \end{abstract}

\maketitle

\noindent \textbf{Keywords.} Brownian motion, capacitance operator, high-contrast resonator system, off-diagonal decay, long-range interactions, nonperiodic structures.\par

\bigskip

\noindent \textbf{AMS Subject classifications.} 31B15, 60J65, 35J05, 35K08, 35P25.

\tableofcontents

\section{Introduction}
Subwavelength resonator systems provide a mathematical framework for the confinement and manipulation of waves on scales smaller than the wavelength. One particular way of realising such subwavelength resonances is through high contrast inclusions embedded in a uniform background medium. In this setting, one can show that, at leading order in the contrast parameter $\delta$, the resonance frequencies and eigenmodes are governed by the eigenvalues and eigenmodes of a generalised capacitance matrix \cite{ammari2018minnaert,ammari2026mathematical}. Minnaert resonances \cite{minnaert1933XVI, ammari2018minnaert} of an air bubble in water are a classical example of such systems. Subwavelength physics also provide a promising research direction for metamaterials design \cite{ammari2026mathematical, davies2025roadmap}.

This connection with metamaterial design motivates the study of large resonator systems. A particularly important question then concerns the off-diagonal decay of the capacitance matrix. Good estimates on the asymptotic decay can be leveraged to obtain results about spectral convergence for finite systems to their infinite counterpart \cite{ammari.davies.ea2025Spectral}, localisation properties of defect modes \cite{ammari.davies.ea2023Convergence} and Anderson localisation \cite{anderson1958Absence, ammari2026long}, where the decay plays a particularly important role \cite{aizenman2015Random}.

Beyond wave physics, the capacitance matrix is a central object in electrostatics, where it is computed for integrated circuits by floating random walk solvers \cite{lecoz1992stochastic, yu2013rwcap}. Here, decay estimates indicate how much computational effort is needed to accurately resolve distant resonator interactions, and when far-away resonators can be neglected.

Results on this decay, however, are scarce. While it is known that the decay is exponential for crystals\cite{ammari2025analysis}, only recently an upper bound was establish for a periodic chain whose fundamental cell contain one resonator \cite{ammari2026long}. The proof of this revolves around analysing the symbol of the associated Laurent operator, hence the approach is restricted to periodic systems only. To the best of our knowledge, no decay estimates are available for disordered configurations or screens.

In this paper, we establish upper and lower bounds for the
off-diagonal entries of the capacitance matrix for three
classes of infinite, not necessarily periodic, configurations
of spherical resonators: chains contained in a cylinder,
screens contained in a slab, and relatively dense
three-dimensional configurations.
For two distinct resonators with centres $x_i$ and $x_j$,
write $d_{ij}=|x_i-x_j|$.
For chains and screens, respectively, we prove
\[
    |\mathcal{C}_{ij}|
    \asymp\frac{1}{d_{ij}\log^2(2+d_{ij})},
    \qquad
    |\mathcal{C}_{ij}|\asymp d_{ij}^{-3},
\]
where $\asymp$ denotes upper and lower bounds by positive
constant multiples. For crystals, we obtain
\[
    C_4e^{-c_4d_{ij}}
    \leq|\mathcal{C}_{ij}|
    \leq C_3e^{-c_3d_{ij}},
\]
with positive constants $C_3,c_3,C_4,c_4$.
The upper and lower exponential rates are allowed to differ.
All constants depend only on the geometric bounds specified
in \cref{thm: ChainDecay,thm: ScreenDecay,thm: CrystalDecay}.
In particular, they are uniform over all distinct index
pairs and all configurations satisfying those fixed bounds.

Our approach uses the classical connection between the
Dirichlet problem and Brownian motion
\cite{morters2010brownian,chung2012brownian}.
The unit potential associated with a resonator is the
probability that Brownian motion first reaches that
resonator before any other. The upper bounds follow from
survival estimates and bounds on the killed transition
density. The lower bounds combine Green's function estimates
in suitable subdomains with path constructions that avoid
the other resonators.


The paper is organised as follows. In Section \ref{sec: Main}, we introduce the setting and state our main results. Section \ref{sec: Preliminaries} collects the necessary facts on Brownian motion and proves the probabilistic representation of the capacitance matrix. In Section \ref{sec: KilledBM}, we study Brownian motion killed on the resonators and relate hitting probabilities to occupation times. Sections \ref{sec: Upper} and \ref{sec: Lower} contain the proofs of the upper and lower bounds, respectively, and we conclude in Section \ref{sec: Conclusion}.

\section{Setting and Main Result}\label{sec: Main}
Let $\mathcal{I}$ be a countable index set and $\{D_i\}_{i\in\mathcal{I}}$ a family of bounded Lipschitz domains in $\R^3$ with pairwise disjoint closures and set $\mathcal{D}=\bigcup_{i\in\mathcal{I}}D_i$. We further write $\Omega$ for $\mathbb{R}^3 \setminus \overline{\mathcal{D}}$. Inside each resonator $D_i$, let the material density and bulk modulus be given by $\rho_i$ and $\kappa_i$ respectively. Outside the resonators, we denote these quantities by $\rho$ and $\kappa$. The interior and exterior wavespeeds are then given by
\begin{equation}
    v_i = \sqrt{\frac{\kappa_i}{\rho_i}}\quad\text{and}\quad v=\sqrt{\frac{\kappa}{\rho}}.
\end{equation}
We assume that there exist $v_{\max}\geq v_{\min}>0$ such that $v_{\max}\geq v_i,v\geq v_{\min}$. If we define the contrast parameter $\delta_i=\frac{\rho_i}{\rho}$ the resonant frequencies $\omega(\delta)\in\C$ are values such that the system below has a nontrivial solution $u$.
\begin{equation}\label{eq: HelmholtzProblem}
    \begin{cases}
        \Delta u + \dfrac{\omega^2}{v^2}\, u = 0
        & \text{in } \Omega, \\[0.8em]
        \Delta u + \dfrac{\omega^2}{v_i^2}\, u
         = 0
        & \text{in } D_i, \quad i \in \mathcal{I}, \\[0.8em]
        u\big|_+ = u\big|_-
        & \text{on } \partial D_i, \quad i \in \mathcal{I}, \\[0.8em]
        \delta_i \dfrac{\partial u}{\partial \nu}\bigg|_+
        = \dfrac{\partial u}{\partial \nu}\bigg|_-
        & \text{on } \partial D_i, \quad i \in \mathcal{I}, \\[0.8em]
        u \text{ is outgoing}.
    \end{cases}
\end{equation}
Here, $\nu$ denotes the unit normal on $\partial D_i$ pointing out of $D_i$ and the signs $+$ and $-$ denote exterior and interior traces,
respectively.
In the interest of readability, we will restrict ourselves to the case $\delta_i\equiv\delta$, however all the main results should be adaptable to the more general case without much difficulty.
In this manuscript, we will be concerned with the \emph{subwavelength} setting where $\delta\rightarrow 0$. To this end, we define \emph{subwavelength resonance frequencies} as follows:
\begin{definition}
    Given $\delta>0$, a subwavelength resonance frequency $\omega(\delta) \in \C$ is such that
    \begin{itemize}
        \item \eqref{eq: HelmholtzProblem} has a nontrivial solution $u$, 
        \item $\omega(\delta)$ is continuous and satisfies $\omega(\delta)=\mathcal{O}(\sqrt{\delta})$.
    \end{itemize}
\end{definition}
This regime is physically interesting, because the wavelength of the associated eigenmodes can be magnitudes larger than the resonators themselves. The perhaps most classical example of this are so-called Minnaert Resonances, which are the resonant frequencies of air bubbles in water \cite{minnaert1933XVI, ammari2018minnaert}. This peculiar property has made the subwavelength regime an attractive platform for metamaterial design \cite{cheben.halir.ea2018Subwavelength, cummer.christensen.ea2016Controlling}. 

In studying the subwavelength setting, the capacitance matrix approach has proven to be very fruitful, for a thorough introduction we refer the reader to \cite{ammari2026mathematical}. To introduce the capacitance matrix, we define $V_i$ to be the solution to
\begin{equation}\label{eq: DirichletProblem}
    \begin{cases}
        \Delta V_i=0,& \text{in } \Omega,\\
        V_i=\delta_{ij} & \text{on } \partial D_j,\\
        V_i(x)=\mathcal{O}(\lvert x\rvert^{-1}).
    \end{cases}
\end{equation}
The Maximum Principle  implies that $V_i$ is unique if it exists.
\begin{remark}
    In previous literature, the existence of $V_i$ has at times been a delicate issue. In the seminal work \cite{ammari2018minnaert}, layer potential techniques were used. This approach however only works when $\lvert \mathcal{I}\rvert <\infty$, as otherwise the single layer potential may fail to be invertible \cite{ammari2026mathematical}. In more recent literature \cite{ammari2025analysis,ammari2026perturbative}, it was defined through the Dirichlet-to-Neumann operator $\mathcal{T}^{\omega}$. 
\end{remark}
One now defines the capacitance matrix as
\begin{equation}\label{eq: CapacitanceFlux}
    \mathcal{C}_{ij}=-\int_{\partial D_j}\partial_\nu V_i\,\mathrm{d}\sigma,\qquad i,j\in\mathcal{I}.
\end{equation}
We will further need the coefficient matrix $\mathcal{V}=\operatorname{diag}(\frac{v_i^2}{\lvert D_i\rvert})$, where $\lvert \cdot \rvert$ denotes the usual Lebesgue measure.
\begin{theorem}[Fundamental Theorem of Subwavelength Physics]\label{thm:fundamental}
    Assume that $\lvert \mathcal{I}\rvert<\infty$. Then the subwavelength resonance problem  \eqref{eq: HelmholtzProblem} admits $\lvert \mathcal{I}\rvert$ resonant frequencies (with positive real parts) that satisfy
    \begin{equation}
        \omega_n = \sqrt{\delta\lambda_n} + \mathcal{O}(\delta),
    \end{equation}
    with associated eigenmode
    \begin{equation}
        u_n(x)=\sum_{i\in\mathcal{I}}(\mathbf{u_n})_iV_i(x) + \mathcal{O}(\sqrt{\delta}),
    \end{equation}
    where $(\lambda_n,\mathbf{u_n})$ is the $n$th eigenvalue-eigenvector pair of $\mathcal{V}\mathcal{C}$.
\end{theorem}

\begin{remark}
    If $\lvert \mathcal{I} \rvert=\infty$, the situation is more complicated. If the underlying resonator configuration is periodic, one approach is to Floquet-Bloch transform Equation \eqref{eq: HelmholtzProblem} directly and study the resulting problem using quasiperiodic Green's functions. This strategy is of course limited to the periodic setting. While one can define the capacitance matrix in the disordered setting (in the crystal case by inverting the Laplacian \cite{ammari2025analysis,ammari2026perturbative} or alternatively by using the coercivity of the Laplacian in Beppo-Levi spaces), a rigorous reduction will be proven in a future paper.
\end{remark}

Our main results will concern the off-diagonal decay of the capacitance matrix. Estimates on this decay have direct implications for the decay of defect modes, as well as spectral convergence of finite systems.

In what follows we will assume $D_i=B_{r_i}(x_i)$, $\mathcal{D}=\bigcup_{i\in\mathcal{I}}D_i$, where $\{r_i\}_{i\in\mathcal{I}}\subset\mathbb{R}_+$ and $\{x_i\}_{i\in\mathcal{I}}\subset\mathbb{R}^3$ satisfy
\begin{equation}\label{eq: Radial Bounds}
    0< r_{\min}\leq r_i\leq r_{\max}<\infty,
\end{equation}
for some $r_{\max}\geq r_{\min}>0$ and, for some constant $d_{\min}$,
\begin{equation}\label{eq: MinimalDistance}
    0<d_{\min}\leq\inf_{i\neq j}\operatorname{dist}(D_i,D_j). 
\end{equation}
Throughout, we write $z=(z_\perp,z_\parallel)\in\R^2\times\R$ for points of $\R^3$ and, likewise, $B=(B_\perp,B_\parallel)$ for the Brownian motion of \cref{sec: BrownianSetup}.

\begin{figure}[H]
\centering
\begin{tikzpicture}[x=0.98cm,y=0.78cm,>=stealth]
    \draw[gray!65] (-6,-2)--(6,-2) (-6,2)--(6,2);
    \draw[dashed,gray!80] (-6,-1)--(6,-1) (-6,1)--(6,1);
    \draw[->,gray] (-6,0)--(6.3,0) node[right] {$z_\parallel$};
    \foreach \x/\y/\r in {-4.5/0.32/0.28,-1.6/-0.28/0.36,1.3/0.36/0.30,4.3/-0.18/0.40}{
        \draw[blue,fill=blue!6] (\x,\y) circle[radius=\r];
    }
    \fill[blue] (-4.5,0.32) circle[radius=0.025];
    \fill[blue] (4.3,-0.18) circle[radius=0.025];
    \draw[blue] (-1.6,-0.28) circle[radius=0.60];
    \node[blue] at (-1.6,-0.28) {$D_i$};
    \node[blue] at (1.3,0.36) {$D_j$};
    \node[above left,blue] at (-2.02,0.14) {$A_i$};
    \node[fill=white,inner sep=2pt] at (4.3,1) {$\partial Z_0$};
    \node[fill=white,inner sep=2pt] at (4.3,2) {$\partial Z_1$};
    \draw[<->,blue] (-1.6,-1.55)--(1.3,-1.55)
        node[midway,fill=white,inner sep=2pt] {$|s_i-s_j|$};
    \draw[<->] (-5.55,0)--(-5.55,1)
        node[midway,left,fill=white,inner sep=1pt] {$1$};
    \draw[<->] (5.65,0)--(5.65,2)
        node[midway,left,fill=white,inner sep=1pt] {$2$};
\end{tikzpicture}
\caption{A schematic axial view after scaling to $\rho\leq 1/2$ (\cref{lem: Scaling}).}
\label{fig: ChainGeometry}
\end{figure}
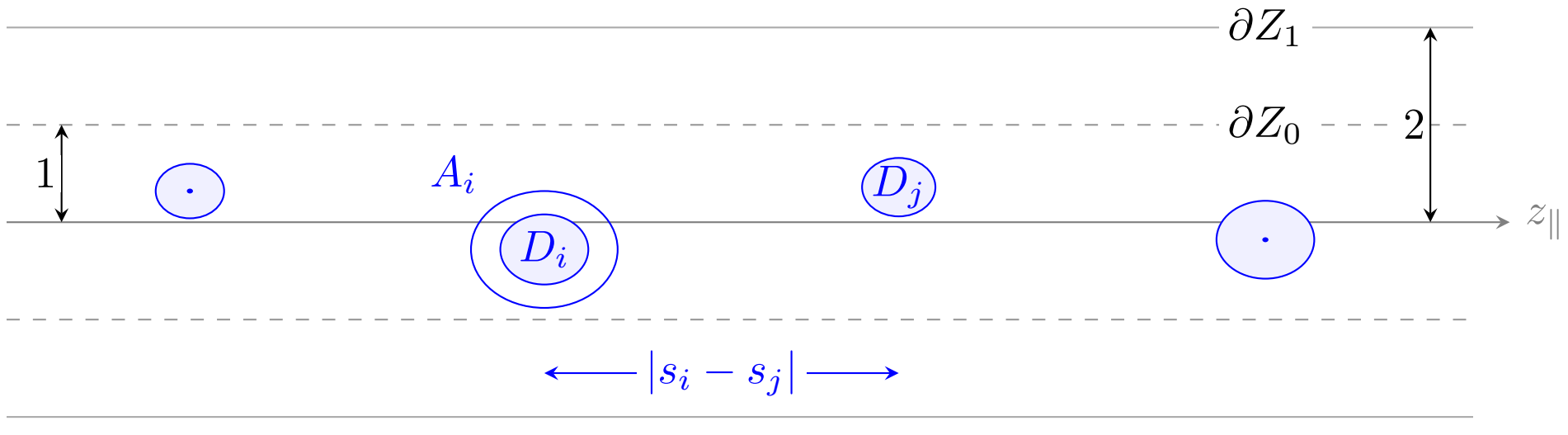

Our first main result concerns chains of spherical resonators. A first result in this direction has recently been proven in \cite{ammari2026long} for a periodic chain whose fundamental cell contains a single resonator. We extend this result to a much larger class of chain-like geometries.
\begin{theorem}[Chain decay]\label{thm: ChainDecay}
    Let $\mathcal{I}=\mathbb{Z}$ and write $x_i = (a_i, s_i)\in\R^2\times\R$ for the centres of the resonators. Assume that additionally to assumptions \eqref{eq: Radial Bounds} and \eqref{eq: MinimalDistance} there exist $L, \rho>0$ such that
    \begin{equation}\label{eq: Cylinder Assumption}
        \lvert a_i\rvert + r_i\leq \rho\quad\text{and}\quad 0\leq s_{i+1} - s_i\leq L.
    \end{equation}
    Then there exist constants $C_1\geq c_1>0$, depending only on $r_{\min}, r_{\max}, d_{\min}, \rho, L$, such that the following bounds hold
    \begin{equation}
        \frac{c_{1}}{d_{ij}\log^2 (2+d_{ij})}\leq \lvert \mathcal{C}_{ij}\rvert\leq\frac{C_{1}}{d_{ij}\log^2 (2+d_{ij})},
    \end{equation}
    for all distinct \(i,j\in \mathcal{I}\), where $d_{ij}=\lvert x_{i}-x_j\rvert$.
\end{theorem}
Equation \eqref{eq: Cylinder Assumption} implies that all resonators are contained in an infinite cylinder of radius $\rho$ along the third coordinate axis. Moreover, as the balls $D_i$ are disjoint with radii at least $r_{\min}$, only finitely many of them meet a bounded part of this cylinder, so that $s_i\to\pm\infty$ as $i\to\pm\infty$; together with $0\leq s_{i+1}-s_i\leq L$ this shows that every point of the axis lies within distance $L$ of some $s_i$. For a visualisation of this setup, see Figure \ref{fig: ChainGeometry}, where $Z_k=\{z:\lvert z_\perp\rvert<2^k\}$ denotes a cylinder whose boundary is given by $\partial Z_k=\{z:\lvert z_\perp\rvert=2^k\}$. 

\begin{figure}[htbp]
\centering
\begin{tikzpicture}[x=0.9cm,y=0.9cm,>=stealth,
    every node/.style={font=\small}]
    \begin{scope}
        \node at (2.25,2.65) {(a) Projection onto $P$};
        \fill[gray!4] (0,-1.4) rectangle (4.5,1.65);
        \draw[gray!40] (0,-1.4) rectangle (4.5,1.65);
        \foreach \x/\y in {0.4/1.15,1.4/1.15,2.6/1.15,
            0.65/0,1.55/0,2.8/0,4.05/0,
            0.35/-1.05,2.15/-1.05,3.35/-1.05,4.15/-1.05}{
            \fill[gray!75] (\x,\y) circle[radius=0.038];
        }
        \fill[blue!60!black] (1.0,-0.8) circle[radius=0.05];
        \node[below right,blue!60!black,inner sep=2pt]
            at (1.0,-0.8) {$a_i$};
        \fill[blue!60!black] (3.6,1.15) circle[radius=0.05];
        \node[above right,blue!60!black,inner sep=2pt]
            at (3.6,1.15) {$a_j$};
        \draw[dashed,blue!55] (2.3,0.35) circle[radius=0.78];
        \draw[blue!60!black] (2.3,0.35)--(2.8,0)
            node[midway,below left,fill=white,inner sep=1pt] {$\le L$};
        \draw[blue!60!black] (2.24,0.29)--(2.36,0.41)
                                   (2.24,0.41)--(2.36,0.29);
        \node[above left,blue!60!black,inner sep=2pt]
            at (2.3,0.35) {$w$};
        \node at (2.25,-1.85) {$a_k\in\R^2$};
    \end{scope}
    \begin{scope}[shift={(5.7,0)}]
        \node at (4.0,2.65) {(b) Normal view of the slab};
        \fill[gray!5] (0,-1.2) rectangle (7.9,1.2);
        \draw[gray!65] (0,-1.2)--(7.9,-1.2)
                       (0,1.2)--(7.9,1.2);
        \draw[dashed,gray!65] (0,0)--(7.9,0);
        \node[fill=white,inner sep=1pt] at (7.75,0) {$P$};
        \foreach \x/\y/\r in {0.65/0.48/0.30,1.0/-0.65/0.28,
            2.95/0.48/0.32,3.6/-0.6/0.29,4.4/0.2/0.41,
            5.0/-0.57/0.27,6.65/-0.7/0.26,7.4/0.25/0.30}{
            \draw[gray!70,fill=gray!15] (\x,\y) circle[radius=\r];
            \fill[gray!75] (\x,\y) circle[radius=0.024];
        }
        \draw[blue!60!black,fill=blue!6] (2.0,-0.42) circle[radius=0.35];
        \draw[blue!60!black] (2.0,-0.42) circle[radius=0.58];
        \node[blue!60!black] at (2.0,-0.42) {$D_i$};
        \node[blue!60!black,fill=white,inner sep=1pt]
            at (1.6,0.25) {$A_i^\delta$};
        \draw[blue!60!black,fill=blue!6] (6.1,0.36) circle[radius=0.42];
        \draw[blue!60!black] (6.1,0.36) circle[radius=0.65];
        \node[blue!60!black] at (6.1,0.36) {$D_j$};
        \node[blue!60!black,fill=white,inner sep=1pt]
            at (5.55,0.95) {$A_j^\delta$};
        \draw[<->,gray!80] (-0.35,-1.2)--(-0.35,1.2)
            node[midway,left,fill=white,inner sep=1pt] {$2\rho$};
        \draw[->,gray!75] (8.05,-1.35)--(8.05,2.25)
            node[above] {$z_\parallel$};
        \node[right,inner sep=2pt] at (8.05,1.2) {$\rho$};
        \node[right,inner sep=2pt] at (8.05,-1.2) {$-\rho$};
        \draw[dotted,blue!50] (2.0,1.2)--(2.0,1.8)
                              (6.1,1.2)--(6.1,1.8);
        \draw[dashed,blue] (0,1.6)--(7.9,1.6);
        \node[right,inner sep=2pt] at (8.05,1.6) {$2$};
        \draw[blue] (2.0,2.0) circle[radius=0.2];
        \draw[blue] (6.1,2.0) circle[radius=0.2];
        \fill[blue!60!black] (2.0,2.0) circle[radius=0.04];
        \fill[blue!60!black] (6.1,2.0) circle[radius=0.04];
        \node[right,blue!60!black] at (2.2,2.0) {$z_i$};
        \node[right,blue!60!black] at (6.3,2.0) {$z_j$};
        \node at (4.05,2.0) {$\mathcal H=\{z_\parallel>2\}$};
        \draw[dashed,gray!60] (2.0,-1.2)--(2.0,-1.75)
                              (6.1,-1.2)--(6.1,-1.75);
        \draw[<->,blue!60!black] (2.0,-1.65)--(6.1,-1.65)
            node[midway,fill=white,inner sep=2pt] {$|a_i-a_j|$};
    \end{scope}
\end{tikzpicture}
\caption{Schematic geometry of a screen of finite thickness.
Left: a portion of the projected centre set.
Right: representative resonators lie in $|z_\parallel|\le\rho$,
with protective shells about $D_i,D_j$. The points
$z_k=x_k+3n$, $n=(0,0,1)$, are the centres of the balls $Q_k$ of
\cref{lem: PotentialLower}, which lie in the halfspace $\mathcal H=\{z_\parallel>2\}$
used for the lower bound (after scaling to $\rho\le1/2$; not to scale).}
\label{fig: ScreenGeometry}
\end{figure}
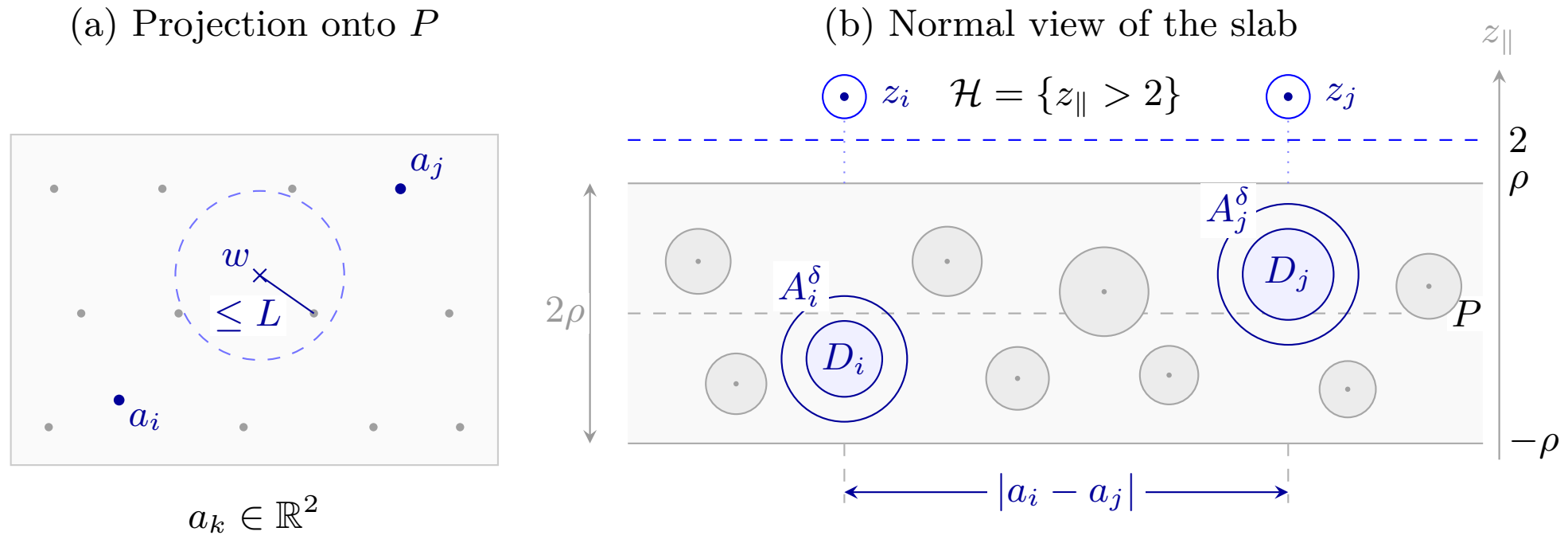

Our second main result concerns screens of resonators. To the best of our knowledge, no prior results exist for this setup, which is illustrated in Figure \ref{fig: ScreenGeometry}.
\begin{theorem}[Screen decay]\label{thm: ScreenDecay}
    Let $\mathcal{I}=\Z$ and again write $x_i = (a_i, s_i)\in\R^2\times\R$ for the centres of the resonators. Assume that additionally to assumptions \eqref{eq: Radial Bounds} and \eqref{eq: MinimalDistance} there exist $L, \rho>0$ such that
    \begin{equation}\label{eq: Slab Assumption}
        \lvert s_i\rvert + r_i\leq \rho\quad\text{and}\quad \sup_{w\in\R^2}\inf_{i\in\mathcal{I}}\lvert a_i-w\rvert\leq L.
    \end{equation}
    Then there exist constants $C_2\geq c_2>0$, depending only on $r_{\min}, r_{\max}, d_{\min}, \rho, L$, such that the following bounds hold
    \begin{equation}
        c_2d_{ij}^{-3}\leq \lvert \mathcal{C}_{ij}\rvert\leq C_2d_{ij}^{-3},
    \end{equation}
    for all distinct \(i,j\in \mathcal{I}\), where again $d_{ij}=\lvert x_i-x_j\rvert$.
\end{theorem}
Equation \eqref{eq: Slab Assumption} implies that all resonators are contained in a slab of width $2\rho$ around the plane $P=\{z:z_\parallel=0\}$ and that every point of $P$ lies within distance $L$ of some projected centre $a_i$, see figure \ref{fig: ChainGeometry} for illustration. 

Finally, we are going to provide a novel proof of exponential decay for crystals of resonators, based on the same strategies employed in the proofs of Theorem \ref{thm: ChainDecay} and \ref{thm: ScreenDecay}, therefore showing that our new probabilistic approach presents a unified framework, not limited to certain geometries.

\begin{theorem}[Crystal decay]\label{thm: CrystalDecay}
    Let $\mathcal{I}=\Z$. Assume that additionally to assumptions \eqref{eq: Radial Bounds} and \eqref{eq: MinimalDistance} there exists $L>0$ such that
    \begin{equation}\label{eq: Crystal Assumption}
        \sup_{z \in \R^3}\inf_{i\in\mathcal{I}}\lvert x_i-z\rvert\leq L.
    \end{equation}
    Then there exist constants $C_3,c_3,C_4,c_4>0$, depending only on $r_{\min}, r_{\max}, d_{\min}, L$, such that the following bound holds
    \begin{equation}
        C_4\exp{\left(-c_4d_{ij}\right)}\leq\lvert \mathcal{C}_{ij}\rvert\leq C_3\exp{\left(-c_3d_{ij}\right)},
    \end{equation}
    where $d_{ij}=\lvert x_i-x_j\rvert$.
\end{theorem}


In the proofs of Theorems \ref{thm: ChainDecay} and \ref{thm: ScreenDecay}, we will make use of the following scaling property, which allows us to restrict ourselves to the case $\rho\leq 1/2$, that is, to configurations contained in the cylinder $Z_0$, respectively the slab $\{z:\lvert z_\parallel\rvert<1\}$.
\begin{lemma}[Scaling]\label{lem: Scaling}
    Let $\lambda>0$ and denote by $V_i^\lambda$ and $\mathcal{C}^\lambda$ the potentials and the capacitance matrix of the rescaled configuration $D_i^\lambda=B_{\lambda r_i}(\lambda x_i)$, $i\in\mathcal{I}$. Then
    \begin{equation}\label{eq: Scaling}
        V_i^\lambda(\lambda z)=V_i(z)\quad\text{for } z\in\Omega\qquad\text{and}\qquad\mathcal{C}^\lambda_{ij}=\lambda\,\mathcal{C}_{ij}\quad\text{for } i,j\in\mathcal{I}.
    \end{equation}
    In particular, it suffices to prove Theorems \ref{thm: ChainDecay} and \ref{thm: ScreenDecay} under the additional assumption $\rho\leq 1/2$.
\end{lemma}
\begin{proof}
    The dilated configuration satisfies \eqref{eq: Radial Bounds} and \eqref{eq: MinimalDistance} with $\lambda r_{\min},\lambda r_{\max},\lambda d_{\min}$ in place of $r_{\min},r_{\max},d_{\min}$, and its exterior domain is $\lambda\Omega$. The function $z\mapsto V_i^\lambda(\lambda z)$ is harmonic in $\Omega$, equals $\delta_{ij}$ on $\partial D_j$ and is $\mathcal{O}(\lvert z\rvert^{-1})$, hence it coincides with $V_i$ by the uniqueness of the solution to \eqref{eq: DirichletProblem}. Differentiating, we obtain $\partial_\nu V_i^\lambda(\lambda y)=\lambda^{-1}\partial_\nu V_i(y)$ for $y\in\partial D_j$. As the surface measure of $\partial D_j^\lambda$ is $\lambda^2$ times that of $\partial D_j$, Equation \eqref{eq: CapacitanceFlux} yields $\mathcal{C}^\lambda_{ij}=\lambda\mathcal{C}_{ij}$.

\end{proof}

\section{Preliminaries}\label{sec: Preliminaries}
\subsection{Brownian motion and stopping times}\label{sec: BrownianSetup}
Throughout, $B$ denotes a standard Brownian motion in $\mathbb{R}^n$. We write $\mathbb{P}_z$ and $\mathbb{E}_z$ for probability and expectation when $B(0)=z$. We use the canonical path space
\[
    \Omega_0=C([0,\infty);\mathbb{R}^n)
\]
and the right-continuous filtration $(\mathcal{F}^+(t))_{t\geq0}$ of \cite[Sections 1.1 and 2.1--2.2]{morters2010brownian}. All stopping times are taken with respect to this filtration. 

For a closed set $H$, we define the stopping time 
\[
    \tau(H)=\inf\{t\geq0:B(t)\in H\},\qquad \inf\varnothing=\infty.
\]
We will often write $\tau:=\tau(\partial\mathcal{D})$, where the intention is clear. Hitting times of closed sets are stopping times, see \cite[Remark 2.14]{morters2010brownian} or \cite[Proposition 1.15]{chung2012brownian}.

We use the strong Markov property in following form \cite[Theorem 1.2]{chung2012brownian}:
\begin{theorem}
    Let $T$ be a stopping time, $E\in\mathcal{F}^+(T)$ and  $f:C([0,\infty);\mathbb{R}^n)\to[0,\infty]$ be measurable. Then it holds that
\begin{equation}\label{eq: StrongMarkov}
    \mathbb{E}_z\left[\mathrm{1}_{E\cap\{T<\infty\}}\,f\left(B(T+\cdot)\right)\right]=\mathbb{E}_z\left[\mathrm{1}_{E\cap\{T<\infty\}}\,\mathbb{E}_{B(T)}[f(B)]\right].
\end{equation}
\end{theorem}

For a measurable functional $F$ of paths we write $F\circ\theta_T=F(B(T+\cdot))$ on $\{T<\infty\}$; for instance, if $T'=\inf\{t\geq0:B(t)\in H\}$, then $T'\circ\theta_T=\inf\{t\geq0:B(T+t)\in H\}$.

\subsubsection{Elementary estimates for Brownian motion}\label{sec: BMEstimates}
We first collect three estimates: a Gaussian displacement bound, an exit-time bound for a disc, and the logarithmic tail of a planar return time.

Let $\overline{\Phi}(x)=\frac{1}{\sqrt{2\pi}}\int_x^\infty e^{-v^2/2}\,\mathrm{d}v$ denote the standard Gaussian tail.
\begin{lemma}[Gaussian estimates]\label{lem: Displacement}
    \begin{enumerate}[label=(\alph*)]
        \item For $x\geq0$, $\overline{\Phi}(x)\leq\frac12e^{-x^2/2}$.
        \item Let $W$ be a standard Brownian motion in $\mathbb{R}^n$, $n\in\mathbb{N}$. For $s>0$ and $a>0$,
        \[
            \mathbb{P}\Big(\sup_{0\leq u\leq s}\lvert W(u)\rvert\geq a\Big)\leq 2n\exp\left(-\frac{a^2}{2ns}\right).
        \]
    \end{enumerate}
\end{lemma}
\begin{proof}
    (a) As $(x+v)^2\geq x^2+v^2$ for $x,v\geq0$,
    \[
        \overline{\Phi}(x)=\frac{1}{\sqrt{2\pi}}\int_0^\infty e^{-(x+v)^2/2}\,\mathrm{d}v\leq\frac{e^{-x^2/2}}{\sqrt{2\pi}}\int_0^\infty e^{-v^2/2}\,\mathrm{d}v=\frac12e^{-x^2/2}.
    \]

    (b) Let $X$ be a standard linear Brownian motion and $M(s)=\max_{u\leq s}X(u)$. For $0<b'<b$, \cite[Theorem 2.21]{morters2010brownian} gives $\mathbb{P}(M(s)\geq b)\leq\mathbb{P}(M(s)>b')=2\mathbb{P}(X(s)>b')$; letting $b'\uparrow b$ and using (a), $\mathbb{P}(M(s)\geq b)\leq2\overline{\Phi}(b/\sqrt{s})\leq e^{-b^2/(2s)}$. Applying this also to $-X$, $\mathbb{P}(\sup_{u\leq s}\lvert X(u)\rvert\geq b)\leq2e^{-b^2/(2s)}$. If $\max_{u\leq s}\lvert W_k(u)\rvert<a/\sqrt{n}$ for all coordinates $k$, then $\max_{u\leq s}\lvert W(u)\rvert^2\leq\sum_{k}\max_{u\leq s}W_k(u)^2<a^2$; a union bound over the $n$ coordinates gives the claim.
\end{proof}
\begin{lemma}[Exit from discs]\label{lem: DiskExit}
    Let $W$ be a planar Brownian motion, $R>0$ and $T_R=\inf\{t\geq0:\lvert W(t)\rvert\geq R\}$. Then $\mathbb{P}_w(T_R>a)\leq\min\{1,9R^2/a\}$ for all $w\in\mathbb{R}^2$ and $a>0$. In particular $T_R<\infty$ almost surely.
\end{lemma}
\begin{proof}
    Let $h=3R^2$. For every $v\in\mathbb{R}^2$ the first coordinate of $W(h)$ under $\mathbb{P}_v$ has density bounded by $(2\pi h)^{-1/2}$, so $\mathbb{P}_v(\lvert W(h)\rvert<R)\leq2R(2\pi h)^{-1/2}=2(6\pi)^{-1/2}<\frac12$. Let $E_N=\{\lvert W(kh)\rvert<R\text{ for }1\leq k\leq N\}$ and $E_0$ the sure event; $E_{N-1}$ is determined by $W(h),\dots,W((N-1)h)$ and hence lies in $\mathcal{F}^+((N-1)h)$. By \eqref{eq: StrongMarkov} at the time $(N-1)h$ with $E=E_{N-1}$, $\mathbb{P}_w(E_N)=\mathbb{E}_w[\mathrm{1}_{E_{N-1}}\mathbb{P}_{W((N-1)h)}(\lvert W(h)\rvert<R)]\leq\frac12\mathbb{P}_w(E_{N-1})$, so $\mathbb{P}_w(T_R>Nh)\leq\mathbb{P}_w(E_N)\leq2^{-N}$. With $N=\lfloor a/h\rfloor$ and $e^{-y}\leq1/y$,
    \[
        \mathbb{P}_w(T_R>a)\leq2^{-N}\leq2\exp\left(-\frac{a\log 2}{3R^2}\right)\leq\frac{6R^2}{a\log2}\leq\frac{9R^2}{a}.\qedhere
    \]
\end{proof}

\begin{lemma}[Planar excursions]\label{lem: Excursion}
    Let $W$ be a planar Brownian motion started in $w$ with $\lvert w\rvert=2$, and $\kappa=\inf\{t\geq0:\lvert W(t)\rvert\leq1\}$. Then $\mathbb{P}_w(\kappa>a)\leq12/\log(2+a)$ for all $a>0$; in particular $\kappa<\infty$ almost surely.
\end{lemma}
\begin{proof}
    Let $m\geq2$ be an integer and $T=T_{2^m}$ as in \cref{lem: DiskExit}. Since $1<\lvert w\rvert<2^m$, by continuity $\kappa$ and $T$ are the hitting times of the circles of radii $1$ and $2^m$, and $T<\infty$ almost surely. By \cite[Theorem 3.18]{morters2010brownian} (planar case), $\mathbb{P}_w(T<\kappa)=1-\frac{\log2^m-\log2}{\log2^m}=\frac1m$. Since $\lvert W(T)\rvert=2^m\neq1=\lvert W(\kappa)\rvert$ whenever both are finite, $\{T=\kappa\}$ is a null set, and $\{\kappa>a\}\subset\{T<\kappa\}\cup\{T>a\}$ up to this null set. By \cref{lem: DiskExit},
    \[
        \mathbb{P}_w(\kappa>a)\leq\frac1m+\frac{9\cdot4^m}{a}.
    \]
    For $a\geq2^{16}$ take $m=\lfloor\log_2(a)/4\rfloor\geq4$; then $4^m\leq a^{1/2}$ and $m\geq\log_2(a)/4-1\geq\log_2(a)/8$, so $\mathbb{P}_w(\kappa>a)\leq8\log2/\log a+9a^{-1/2}$. Since $v\mapsto\log(v)/\sqrt v$ decreases on $[e^2,\infty)$, $\log(a)/\sqrt a\leq16\log2/2^8<0.044$, hence $9a^{-1/2}\leq0.4/\log a$ and $\mathbb{P}_w(\kappa>a)\leq6/\log a\leq12/\log(2+a)$, using $2+a\leq a^2$. For $0<a<2^{16}$ the bound is trivial, as $12/\log(2+a)>12/\log(2^{17})>1$.
\end{proof}


\subsection{Exterior harmonic extension}
In this subsection, we will prove an adaption of \cite[Theorem 1.23]{chung2012brownian} to exterior domains for functions with suitable decay.

\begin{theorem}[Exterior harmonic extension]\label{thm: HarmonicExtension}
    Let $n\geq3$ and let $\{D_i\}_{i\in\mathcal I}$ be a locally finite, countable family of bounded Lipschitz domains in $\mathbb R^n$ with pairwise disjoint closures. Set $\mathcal D=\bigcup_{i\in\mathcal I}D_i$ and $\Omega=\mathbb R^n\setminus\overline{\mathcal D}$.  Write $\tau=\tau(\partial\mathcal D)$. For any continuous function $\varphi:\partial\mathcal D\to\mathbb R$ satisfying
    \begin{equation}
        \lim_{R\rightarrow\infty}\sup_{x\in \partial \mathcal{D}\cap B_R(0)^c}\lvert \varphi(x)\rvert=0
    \end{equation}
    the function defined by
    \begin{equation}
        u(x)=\mathbb{E}_x\left[\varphi(B(\tau(\partial \mathcal{D})))\mathrm{1}_{\{\tau<\infty\}}\right]
    \end{equation}
    for $x\in\Omega$, extended by $\varphi$ on $\partial\Omega$, is the unique bounded function in $C(\overline\Omega)$ that is harmonic in $\Omega$, has boundary values $\varphi$, and vanishes at infinity.
\end{theorem}
\begin{proof}   
    The continuity of $\varphi$ and its decay at infinity imply that it is bounded. Existence and continuity hence follows from \cite[Theorem 1.23]{chung2012brownian} since the Lipschitz condition ensures regularity of the boundary by the exterior cone criterion \cite[Proposition 1.22]{chung2012brownian}.
    Next, we show that $u$ vanishes at infinity. For any $\mathcal{J}\subset \mathcal{I}$ finite define
    \begin{equation}
        \varphi^\mathcal{J}(x)=\begin{cases}
            \varphi(x),\quad x\in\partial\mathcal{D}\cap\bigcup_{i \in \mathcal{J}}\partial D_i,\\
            0,\quad \text{else}.
        \end{cases}
    \end{equation}
    Fix $\epsilon>0$. By local finiteness and the decay of $\varphi$ at infinity, choose a finite subset $\mathcal{J}\subset \mathcal{I}$ such that $\lvert\varphi-\varphi^\mathcal{J}\rvert<\epsilon$. Define
    \begin{equation}
        u^\mathcal{J}(x)=\mathbb{E}_x\left[\varphi^\mathcal{J}(B(\tau(\partial \mathcal{D})))\mathrm{1}_{\{\tau<\infty\}}\right].
    \end{equation}
    Choose $R_\epsilon$ so that the support of $\varphi^{\mathcal J}$ is contained in $B_{R_\epsilon}(0)$. A path contributing to $u^{\mathcal J}(x)$ must hit this ball. The Brownian hitting probability of a ball in dimension $n\geq3$ therefore gives, for $\lvert x\rvert>R_\epsilon$,
    \[
        \lvert u^{\mathcal J}(x)\rvert
        \leq\lVert\varphi\rVert_\infty
        \left(\frac{R_\epsilon}{\lvert x\rvert}\right)^{n-2}.
    \]
    Since $\lVert u-u^{\mathcal J}\rVert_\infty\leq\epsilon$, it follows that
    \begin{equation}
        \lvert u(x)\rvert\leq\epsilon+
        \lVert\varphi\rVert_\infty
        \left(\frac{R_\epsilon}{\lvert x\rvert}\right)^{n-2}.
    \end{equation}
    Letting $\lvert x\rvert\to\infty$ and then $\epsilon\downarrow0$ proves the decay of $u$.
    
    Finally, we show uniqueness. Let $v$ be another solution with the same boundary values and decay at infinity. Then $h=u-v$ is harmonic in the exterior domain, vanishes on $\partial\mathcal D$, and tends to zero at infinity. Apply the maximum principle to $h$ and $-h$ on each component of
    $(\mathbb{R}^n\setminus\overline{\mathcal D})\cap B_R(0)$, and let $R\to\infty$.
\end{proof}

\begin{lemma}\label{lem: Connectedness}
    Suppose the family $\{D_i\}_{i\in\mathcal{I}}$ satisfies \eqref{eq: Radial Bounds} and \eqref{eq: MinimalDistance}. Then $\{D_i\}_{i\in\mathcal{I}}$ is locally finite and $\Omega$ is connected. In particular, the assumptions of Theorem \ref{thm: HarmonicExtension} are satisfied.
\end{lemma}
\begin{proof}
    Because $r_i\leq r_{\max}$, any $D_i$ intersecting $B_R(0)$ for some $R>0$ must be contained in $B_{R+2r_{\max}}(0)$. However, $\lvert D_i\rvert\geq 4\pi r_{\min}^3/3$ and hence any bounded set can contain only finitely many inclusions. To see that $\Omega$ is connected, cover each $D_i$ with $B_{r_i + d_{\min}/3}(x_i)$. Let now $x,y\in\Omega$ be two arbitrary points and $\gamma$ the straight segment connecting them, which meets only finitely many of these balls. Whenever $\gamma$ enters some $B_{r_i + d_{\min}/3}(x_i)$, we may connect the entry and exit point by an arc on $\partial B_{r_i + d_{\min}/3}(x_i)\subset\Omega$; if $x$ or $y$ lies inside such a ball, we first move it radially to the sphere, which is possible as $x,y\notin\overline{\mathcal{D}}$. This yields a path $\gamma'\subset \Omega$.
\end{proof}
We thus obtain the following Corollary
\begin{corollary}\label{cor: NewtonianPotential}
    Under the assumptions \eqref{eq: Radial Bounds} and \eqref{eq: MinimalDistance}, the unique solution to the problem \eqref{eq: DirichletProblem} is given by
    \begin{equation}\label{eq: NewtonianPotential}
    V_i(x) = \mathbb{E}_x\left[\mathrm{1}_{\tau<\infty}\mathrm{1}_{\partial D_i}(B(\tau))\right]=\mathbb{P}_x(\tau<\infty, B(\tau)\in\partial D_i).
\end{equation}
\end{corollary}
As an immediate consequence, one obtains the following decay estimate.
\begin{proposition}\label{prop: PotentialBound}
    For every $z\in\Omega$,
    \begin{equation}
        0\leq V_i(z)\leq\frac{r_i}{\lvert z- x_i\rvert}.
    \end{equation}
\end{proposition}
\begin{proof}
    Let $T=\inf\{t>0:\lvert B(t)-x_i\rvert=r_i\}$. On $\{\tau<\infty,\ B(\tau)\in\partial D_i\}$ we have $T\leq\tau<\infty$. Since $\lvert z-x_i\rvert>r_i$, \cite[Corollary 3.19]{morters2010brownian} (three-dimensional case), applied to $B-x_i$, gives $V_i(z)\leq\mathbb{P}_z(T<\infty)=r_i/\lvert z-x_i\rvert$.
\end{proof}
\subsection{From flux to hitting probabilities}
The next lemma is the final step in going from estimating a flux to upper bounding hitting probabilities. First however, we need to define the open shells 
\begin{equation}
    A_j^\delta=\{z\in\R^3:r_{j}<\lvert x_{j}-z\rvert<\delta + r_{j}\}.
\end{equation}
Note for $\delta \in (0,d_{\min}/2]$ we have \begin{equation}\label{eq: Shells}
    \begin{gathered}        \overline{A^\delta_i}\cap\overline{D_k}=\emptyset\ (k\neq i),\qquad\overline{A^\delta_i}\cap\overline{D_i}=\partial D_i,\\
        \text{in particular}\quad A^\delta_i\subset\Omega,\quad\partial B_{r_i+\delta}(x_i)\subset\Omega.
    \end{gathered}
\end{equation}

\begin{lemma}[Flux bound]\label{lem: Flux}
    Let $j\in\mathcal{I}$ and $\delta\in(0,d_{\min}/2]$.
    \begin{enumerate}[label=(\alph*)]
        \item Let $u$ be harmonic on $A^\delta_j$ and $M\geq0$ with
        \begin{equation}\label{eq: Barrier}
            0\leq u(z)\leq M\,\frac{r_j^{-1}-\lvert z-x_j\rvert^{-1}}{r_j^{-1}-(r_j+\delta)^{-1}},\qquad z\in A^\delta_j.
        \end{equation}
        Then $\Psi_u(\varrho)=\int_{\partial B_\varrho(x_j)}\partial_ru\,\mathrm{d}\sigma$ does not depend on $\varrho\in(r_j,r_j+\delta)$, and its value $\Psi_u$ satisfies $0\leq\Psi_u\leq4\pi r_j(r_j+\delta)M/\delta$.
        \item For $i\neq j$, $u=V_i$ satisfies \eqref{eq: Barrier} with $M=M^\delta_{ij}=\sup\{V_i(z):\lvert z-x_j\rvert=r_j+\delta\}$, and $u=1-V_j$ satisfies \eqref{eq: Barrier} with $M=1$.
    \end{enumerate}
\end{lemma}
\begin{proof}
    \emph{(a)} For $r_j<\varrho_1<\varrho_2<r_j+\delta$ the divergence theorem on the closed shell $\{\varrho_1\leq\lvert z-x_j\rvert\leq\varrho_2\}\subset A^\delta_j$ gives $\Psi_u(\varrho_2)-\Psi_u(\varrho_1)=\int\Delta u=0$. This establishes \emph{(a)} for $\varrho\in(r_j,r_j+\delta)$. 

    Let $\mu(\varrho)=\frac{1}{4\pi}\int_{\partial B_1(0)}u(x_j+\varrho\omega)\,\mathrm{d}\sigma(\omega)$. Since $\nabla u$ is bounded on closed shells contained in $A^\delta_j$, we may differentiate under the integral and get $\mu'(\varrho)=\frac{1}{4\pi}\int_{\partial B_1(0)}\nabla u(x_j+\varrho\omega)\cdot\omega\,\mathrm{d}\sigma(\omega)=\Psi_u/(4\pi\varrho^2)$, so $\mu(\varrho)=\mu(\varrho_0)+\frac{\Psi_u}{4\pi}(\varrho_0^{-1}-\varrho^{-1})$ for $r_j<\varrho_0<\varrho<r_j+\delta$.

    By \eqref{eq: Barrier}, $0\leq\mu(\varrho_0)\leq M(r_j^{-1}-\varrho_0^{-1})/(r_j^{-1}-(r_j+\delta)^{-1})\to0$ as $\varrho_0\downarrow r_j$, hence
    \begin{equation}\label{eq: SphericalMean}
        \mu(\varrho)=\frac{\Psi_u}{4\pi}\left(\frac{1}{r_j}-\frac1\varrho\right),\qquad\varrho\in(r_j,r_j+\delta).
    \end{equation}
    As $\mu\geq0$, $\Psi_u\geq0$; comparing \eqref{eq: SphericalMean} with the bound on $\mu(\varrho)$ from \eqref{eq: Barrier} gives $\Psi_u\leq4\pi M/(r_j^{-1}-(r_j+\delta)^{-1})=4\pi r_j(r_j+\delta)M/\delta$.
    
    To extend the statement to $\varrho=r_j$, fix $\delta'\in(0,\delta)$ and set $A=\{z:r_j<\lvert z-x_j\rvert<r_j+\delta'\}$. By \eqref{eq: Barrier}, $u\in C^0(\overline{A})\cap C^\infty(A)$ with boundary values $\varphi\in C^\infty(\overline{A})$, where $\varphi$ vanishes near $\partial D_j$ and coincides with $u$ near the outer sphere, so that $u\in C^{2,\alpha}(\overline{A})$ by classical regularity theory \cite[Theorem 6.19]{gilbarg1977elliptic}. Hence $\nabla u$ is continuous to the boundary, and $\Psi_u(r_j)=\lim_{\varrho\downarrow r_j}\Psi_u(\varrho)=\Psi_u$.
    \emph{(b)} Define
    \[
        \displaystyle H_j^\delta(z)
        =\frac{r_j^{-1}-|z-x_j|^{-1}}
               {r_j^{-1}-(r_j+\delta)^{-1}}.
    \]
    This function is harmonic on $A_j^\delta$, equals zero on the inner sphere, and equals one on the outer sphere. By \cref{thm: HarmonicExtension}, both $V_i$ for $i\ne j$ and $1-V_j$ are harmonic on the shell, continuous on its closure, and zero on the inner sphere. On the outer sphere they are bounded above by $M_{ij}^\delta$ and $1$, respectively. The maximum principle applied to $V_i-M_{ij}^\delta H_j^\delta$ and to $1-V_j-H_j^\delta$ proves \eqref{eq: Barrier} in both cases.
\end{proof}
Consequently, applying \cref{lem: Flux}(a) with the constants from \cref{lem: Flux}(b) to Equation \eqref{eq: CapacitanceFlux}, we thus obtain
\begin{equation}\label{eq: CapacitanceByPotential}
    \lvert\mathcal{C}_{ij}\rvert\leq\frac{4\pi r_j(r_j+\delta)}{\delta}\,M^\delta_{ij},\qquad i\neq j,
\end{equation}
thereby successfully going from having to estimate the flux to having to estimate hitting probabilities on a shell around the inclusions. As an immediate application we recover the classical $1/d_{ij}$ decay and can control the effect of removing distant resonators.

\begin{corollary}\label{cor: CrudeDecay}
    For $i\neq j$ and $\delta\in(0,d_{\min}/2]$,
    \begin{equation}
        \lvert\mathcal{C}_{ij}\rvert\leq\frac{4\pi r_ir_j(r_j+\delta)}{\delta}\,\frac{1}{\lvert x_i-x_j\rvert-(r_j+\delta)}.
    \end{equation}
\end{corollary}
\begin{proof}
    Since $\operatorname{dist}(D_i,D_j)=\lvert x_i-x_j\rvert-r_i-r_j$, \eqref{eq: MinimalDistance} gives $\lvert x_i-x_j\rvert\geq r_i+r_j+d_{\min}$, so every $z$ with $\lvert z-x_j\rvert=r_j+\delta$ satisfies $\lvert z-x_i\rvert\geq\lvert x_i-x_j\rvert-r_j-\delta\geq r_i+d_{\min}-\delta>0$. By \cref{prop: PotentialBound}, $M^\delta_{ij}\leq r_i/(\lvert x_i-x_j\rvert-r_j-\delta)$, and \eqref{eq: CapacitanceByPotential} gives the claimed bound.
\end{proof}
\begin{proposition}\label{prop: Truncation}
    Let $\emptyset\neq\mathcal{J}\subset\mathcal{I}$, and denote by $\mathcal{C}^{\mathcal{J}}$ and $\mathcal{C}^{\mathcal{I}}=\mathcal{C}$ the capacitance matrices \eqref{eq: CapacitanceFlux} of $\mathcal{D}^{\mathcal{J}}=\bigcup_{k\in\mathcal{J}}D_k$ and $\mathcal{D}^{\mathcal{I}}=\mathcal{D}$. Let $i,j\in\mathcal{J}$ and $R=\inf_{m \in \mathcal{I}\setminus\mathcal{J}}\operatorname{dist}(x_i,D_m)$ (with $R=\infty$ if $\mathcal{J}=\mathcal{I}$). Then for every $\delta\in(0,d_{\min}/2]$
    \begin{equation}
        0\leq\mathcal{C}^{\mathcal{I}}_{ij}-\mathcal{C}^{\mathcal{J}}_{ij}\leq\frac{4\pi r_ir_j(r_j+\delta)}{\delta R}.
    \end{equation}
\end{proposition}
\begin{proof}
    The case $\mathcal J=\mathcal I$ is immediate. Otherwise, set $h=V_i^{\mathcal J}-V_i^{\mathcal I}$ on $\Omega^{\mathcal I}$. By \cref{thm: HarmonicExtension}, $h$ is harmonic in $\Omega^{\mathcal I}$, continuous on its closure, and vanishes at infinity.

   On a retained sphere $\partial D_k$, $k\in\mathcal J$, both potentials have the same boundary value, so $h=0$. On an omitted sphere $\partial D_m$, $m\notin\mathcal J$, we have $V_i^{\mathcal I}=0$. The one-ball bound of \cref{prop: PotentialBound}, applied to the family $\mathcal J$, gives
    \[
        \displaystyle 0\leq h(w)=V_i^{\mathcal J}(w)
        \leq\frac{r_i}{|w-x_i|}\leq\frac{r_i}{R},
        \qquad w\in\partial D_m.
    \]
    Apply the maximum principle to $h$ and $-h$ on the components of $\Omega^{\mathcal I}\cap B_{\varrho}(0)$, and let $\varrho\to\infty$. Since $h$ vanishes at infinity, this yields $0\leq h\leq r_i/R$ throughout $\Omega^{\mathcal I}$.

     On $A_j^\delta$, the function $h$ is harmonic and continuous up to both spheres, with inner boundary value zero and outer boundary values at most $r_i/R$. Comparing it with the radial barrier $H_j^\delta$ from the proof of \cref{lem: Flux}(b) gives
    \[
        \displaystyle 0\leq h(z)\leq\frac{r_i}{R}H_j^\delta(z),
        \qquad z\in A_j^\delta.
    \]
    Thus \cref{lem: Flux}(a) applies with $M=r_i/R$. By linearity of the flux and \eqref{eq: CapacitanceFlux},
    \[
        \displaystyle 0\leq\mathcal C^{\mathcal I}_{ij}-\mathcal C^{\mathcal J}_{ij}
        =\Psi_h\leq\frac{4\pi r_ir_j(r_j+\delta)}{\delta R}.\qedhere
    \]
\end{proof}

\section{Killed Brownian Motion}\label{sec: KilledBM}
We begin this section by introducing the probability density of killed Brownian motion. This then allows us to define the Green's function for the exterior domain $\Omega$, which will then be used to obtain a quantitative estimate on $M_{ij}$ from \cref{lem: Flux}.
\subsection{The killed transition density and Green's function}\label{sec: KilledDensity}
The killed transition density describes Brownian motion before its first contact with a resonator. We follow \cite[Chapter 2]{chung2012brownian}.

Let $p_t(z)=(2\pi t)^{-3/2}e^{-\lvert z\rvert^2/(2t)}$ denote the Gaussian kernel in $\mathbb{R}^3$ and define
\begin{equation}
    r^\Omega_t(x,y)=\mathbb{E}_x[\mathrm{1}_{\tau<t}p_{t-\tau}(B(\tau) - y)].
\end{equation}
The killed probability density is then defined as
\begin{equation}
    p^\Omega_t(x,y)=p_t(x-y) - r^\Omega_t(x,y).
\end{equation}
The following result can be found in \cite{hunt1956some}, \cite[Theorem 2.4]{chung2012brownian}.

\begin{theorem}\label{thm: KilledDensity}
    Let $\Omega$ be a domain in $\R^d$. Then for any $t>0$, $x\in\Omega$ and $A$ in the Borel $\sigma$-algebra of $\R^d$
    \begin{equation}
        \mathbb{P}_x(\tau>t, B(t)\in A)=\int_Ap^\Omega_t(x,y)dy.
    \end{equation}
    The function $p^\Omega_t$ is symmetric on $\R^d\times\R^d$, continuous on $(\R^d\setminus\partial\Omega)\times (\R^d\setminus\partial\Omega)$, and strictly positive on $\Omega \times \Omega$.
    For any $t>s>0,x,y\in\mathbb{R}^d$, we have
    \begin{equation}
        p_t^\Omega(x,y)=\int_{\R^d}p_s^\Omega(x,z)p_{t-s}^\Omega(z,y)dz.
    \end{equation}
    For any $t>0,y \in\Omega$ and regular point $z\in\partial\Omega$, we have
    \begin{equation}
        \lim_{\substack{x\rightarrow z\\x\in\Omega}}p_t^\Omega(x,y)=0.
    \end{equation}
\end{theorem}

With the killed transition density in hand, we can now define the killed Green's function and operator
\begin{align}    
    \mathcal{G}^\Omega(x,y) &= \int_0^\infty p_t^\Omega(x,y)dt\\
    \mathcal{G}^\Omega f(x)&=\int_{\R^d} \mathcal{G}^\Omega (x,y)f(y)dy = \mathbb{E}_x\left[\int_0^\tau f(B(t))dt\right].
\end{align}
Finally, writing $\mathcal{G}^0(x,y)=\frac{1}{2\pi\lvert x-y\rvert}$ for the free Green's function, the following identity holds \cite[Theorem 2.5]{chung2012brownian}:
\begin{equation}
    \mathcal{G}^\Omega(x,y) = \mathcal{G}^0(x,y)-\mathbb{E}_x\left[\mathrm{1}_{\tau<\infty}\mathcal{G}^0(B(\tau),y)\right].
\end{equation}
By an abuse of notation, we define 
\begin{equation}\label{eq: OccupationTime}
        \mathcal{G}^\Omega(x,\Lambda)=\mathcal{G}^\Omega \mathrm{1}_\Lambda(x)
        = \mathbb{E}_x\left[\int_0^\tau \mathrm{1}_\Lambda(B(t))dt\right]= \int_0^\infty\mathbb{P}_x\left(\tau>t,B(t)\in\Lambda\right)dt,
\end{equation}
where the last equality follows from Theorem \ref{thm: KilledDensity}. It is not difficult to see that $\mathcal{G}^\Omega(x,\Lambda)$ is thus the expected occupation time in $\Lambda$ of a Brownian motion started in $x$ until it is killed. 

In \cref{lem:occupation}, we will use the occupation time $\mathcal{G}^\Omega(\cdot,A_i^\delta)$ to upper bound the potential $V_i$. Consequently, we require good control on $\mathcal{G}^\Omega(\cdot,A_i^\delta)$, which will be done through estimating the killed kernel density. This is the first time where geometry will enter the picture. To avoid essentially proving the same result twice, we state a general version, which will later be adapted to the respective setup.

\begin{lemma}\label{lem:split}
Let $\Pi$ be an orthogonal projection and $x,y \in \Omega$ with $R=|\Pi x-\Pi y|>0$. Then it holds that

\begin{equation}\label{eq:split}
    \begin{aligned}
        p_t^\Omega(x,y)&\le (2 / \pi t)^{3/2}
        [\mathbb{P}_y(\tau>t/4)\mathbb{P}_x\{\tau>t/2,\ |\Pi(B(t/2)-x)|\ge R/2\}\\
        &+\mathbb{P}_x(\tau>t/4)\mathbb{P}_y\{\tau>t/2,\ |\Pi(B(t/2)-y)|\ge R/2\}].
    \end{aligned}
\end{equation}
For any $p>0$ and $\eta>0$,
\begin{equation}\label{eq:moment-tail}
 \mathbb{P}_z\left(\tau>s,\ |\Pi(B(s)-z)|\ge \eta\right)\le\min\{\mathbb{P}_z(\tau>s),\ \eta^{-p}\mathbb{E}_z|\Pi(B(\tau\wedge s)-z)|^p\}.
\end{equation}
\end{lemma}
\begin{proof}
Let $E_x=\{u:|\Pi(u-x)|\ge R/2\}$ and define $E_y$ similarly.  Since $|\Pi x-\Pi y|=R$, the triangle inequality implies $E_x\cup E_y=\R^3$. \cref{thm: KilledDensity} now implies

\begin{equation}\label{eq:split1}
    p_t^\Omega(x,y)\leq\int_{E_x}p_{t/2}^\Omega(x,z)p_{t/2}^\Omega(z,y)dz + \int_{E_y}p_{t/2}^\Omega(x,z)p_{t/2}^\Omega(z,y)dz.
\end{equation}

Applying the semigroup property once more, we have

\begin{equation}\label{eq:split2}
    \begin{aligned}
        p_{t/2}^\Omega(x,y)&=\int_{\mathbb{R}^3}p_{t/4}^\Omega(x,z)p_{t/4}^\Omega(z,y)dz \\
        &\leq\left(2/\pi t\right)^{3/2}\mathbb{P}_{y}(\tau>t/4).
    \end{aligned}
\end{equation}
Combining the above equation with Equation \eqref{eq:split1}, we obtain the Equation \eqref{eq:split}:
\begin{equation}
    \begin{aligned}
        p_t^\Omega(x,y)&\le (2 / \pi t)^{3/2}
        [\mathbb{P}_y(\tau>t/4)\int_{E_x}p_{t/2}^\Omega(x,z)dz +\mathbb{P}_x(\tau>t/4)\int_{E_y}p_{t/2}^\Omega(y,z)dz]\\
        &=(2 / \pi t)^{3/2}
        [\mathbb{P}_y(\tau>t/4)\mathbb{P}_x\{\tau>t/2,\ |\Pi(B(t/2)-x)|\ge R/2\}\\
        &+\mathbb{P}_x(\tau>t/4)\mathbb{P}_y\{\tau>t/2,\ |\Pi(B(t/2)-y)|\ge R/2\}].
    \end{aligned}
\end{equation}
It thus remains to prove Equation \eqref{eq:moment-tail}. The estimate
\begin{equation}
     \mathbb{P}_z\left(\tau>s,\ |\Pi(B(s)-z)|\ge \eta\right)\le\mathbb{P}_z(\tau>s),
\end{equation}
is trivial. We thus concern ourselves with the second bound. To this end, note that
\begin{equation}
    \begin{aligned}
        \mathbb{P}_z\left(\tau>s,\ |\Pi(B(s)-z)|\ge \eta\right)&=\mathbb{P}_z\left(\tau>s,\ |\Pi(B(\tau\wedge s)-z)|\ge \eta\right)\\
        &\leq \mathbb{P}_z\left(|\Pi(B(\tau \wedge s)-z)|\ge \eta\right)\\
        &\leq \mathbb{E}_z\left[|\Pi(B(\tau \wedge s)-z)|^p\right]\eta^{-p},
    \end{aligned}
\end{equation}
where we used a Chebyshev inequality in the last line.
\end{proof}

\subsection{From hitting probability to occupation time}

\begin{lemma}\label{lem:occupation}
For $A_i^\delta$ with $\delta \in (0,d_{\min}/2]$ and $z\in\Omega$ satisfying $ |z-x_i|\ge r_i+\delta/2$ it holds that
\begin{equation}\label{eq:occupation}
 V_i(z)\leq {12\over\delta^2}\mathcal G^\Omega(z,A_i^\delta),
 \qquad |z-x_i|\ge r_i+\delta/2.
\end{equation}
\end{lemma}
\begin{proof}
Let $T$ be the first hit of the sphere of radius $r_i+\delta/2$
about $x_i$. Continuity implies that every path first absorbed at
$D_i$ has $T<\tau$, including $T=0$ if the starting point is on
this sphere.

If $w$ is on the sphere of radius $r_i+\delta/2$, then $B_{\delta/2}(w)$ lies
in $A_i^\delta$. Write $\sigma$ for the exit time of $B$ from $B_{\delta/2}(w)$. Optional
stopping of $|B(t)-w|^2-3t$ at $\sigma$ gives
$3\mathbb{E}_w(\sigma)=\mathbb{E}_w|B(\sigma)-w|^2=\delta^2/4$. Hence $\mathbb{E}_w\sigma=\delta^2/12$.
Before $\sigma$ no absorption occurs, and the whole time is spent
in $A_i^\delta$. Strong Markov at $T$, on the event
$\{T<\tau\}\in\mathcal F^+(T)$, therefore gives
\begin{equation}
    \begin{aligned}
        \mathcal G^\Omega(z,A_i^\delta) &=\mathbb{E}_z\left[\int_0^\tau\mathrm{1}_{A_i^\delta}(B(t))dt\right]\\
        &\geq \mathbb{E}_z\left[\mathrm{1}_{T<\tau}\int_T^\tau\mathrm{1}_{A_i^\delta}(B(t))dt\right]\\
        &= \mathbb{E}_z\left[\mathrm{1}_{T<\tau}\mathbb{E}_{B(T)}\left[\int_0^\tau\mathrm{1}_{A_i^\delta}(B(t))dt\right]\right]\\
        &\geq \frac{\delta^2}{12}\mathbb{P}_z\left(T<\tau\right),
    \end{aligned}
\end{equation}
which yields the claim. 
\end{proof}

\section{Upper bounds}\label{sec: Upper}
In this section, we will derive the upper bounds of \cref{thm: ChainDecay,thm: ScreenDecay,thm: CrystalDecay}. While the arguments are morally the same, the estimates have to be adapted to the respective geometry, hence the split into subsections. By \cref{lem: Scaling}, we may assume $\rho\leq 1/2$ throughout this and the following section. We denote by $\Pi_\parallel z=z_\parallel$ and $\Pi_\perp z=z_\perp$ the orthogonal projections onto the chain axis and onto the screen plane, respectively. In each subsection, $q_-$ denotes the constant of the corresponding killing lemma (Lemmas \ref{lem: Killing}, \ref{lem: KillingScreen} and \ref{lem: KillingCrystal}). Roughly, the argument is as follows
\begin{enumerate}[label=\arabic*.]
    \item Derive a killing probability that holds uniformly for points on the boundary of a cylinder/slab.
    \item Use this uniform killing to prove that the number of visits to a fixed shell before absorption has a geometric tail.
    \item Turn the geometric tail of the visits into a tail estimate for $\mathbb{P}_z(\tau>t)$.
    \item Use the tail estimate for $\mathbb{P}_z(\tau>t)$ to obtain an estimate for the killed density $p_t^\Omega$.
    \item Integrate in time to obtain an upper bound for the Green's function $\mathcal{G}^\Omega(x,y)$.
\end{enumerate}
\subsection{Chains}
As outlined above, we first derive a uniform lower bound on the probability of absorption during a visit to $\partial Z_0$.
\begin{lemma}\label{lem: Killing}
    Let $\alpha=\inf\{t\geq0:\lvert B_\perp(t)\rvert\geq2\}$ be the exit time from $Z_1$. For every $z\in\partial Z_0$,
    \begin{equation}\label{eq: qminus}
        \mathbb{P}_z(\tau<\alpha)\geq q_-\in(0,1).
    \end{equation}
    The constant $q_-$ depends only on $r_{\min}$ and $L$.
\end{lemma}
\begin{proof}
    Fix $T=1$ and $\epsilon=\frac{1}{2}\min\{r_{\min},1\}$. The assumptions of Theorem \ref{thm: ChainDecay} imply that for any $z \in\partial Z_0$, there exists some $x_i$ such that $d=\lvert x_i-z\rvert\leq\sqrt{L^2 + 9/4}$. Set $F(t) = z + \frac{t(x_i-z)}{T}$ and consider the event $E=\{\max_{t\leq T}\lvert B(t) - F(t)\rvert<\epsilon \}$. That is, $B$ stays within a tube of radius $\epsilon$ of the path $F$ until time $T$. It is not difficult to see that
    \begin{equation}
        \{\tau<\alpha\}\supset E,
    \end{equation}
    since on $E$ we have $\lvert B(T)-x_i\rvert<\epsilon<r_i$, i.e.\ $B(T)\in D_i$, while $\lvert B_\perp(t)\rvert\leq\max\{\lvert z_\perp\rvert,\lvert a_i\rvert\}+\epsilon<2$ for all $t\leq T$. Hence $\mathbb{P}_z(\tau<\alpha)\geq \mathbb{P}_z(E)$. The Cameron-Martin Theorem \cite[Theorem 1.38, Remark 1.43]{morters2010brownian} yields
    \begin{equation}
        \mathbb{P}_z(E)=\mathbb{E}_0\left[\mathrm{1}_A\exp{\left(\int_0^T F'dB - \frac{1}{2}\int_0^T\lvert F'\rvert^2(s)ds\right)}\right],
    \end{equation}
    where $A=\{\max_{t\leq T}\lvert B(t)\rvert<\epsilon\}$. Noting that $A$ is symmetric with respect to negation, i.e. $B\in A\iff -B\in A$, we obtain
    \begin{equation}
        \begin{aligned}
                \mathbb{P}_z(E)&=\mathbb{E}_0\left[\mathrm{1}_A\cosh\left(\frac{\langle x_i-z,B(T)\rangle}{T}\right)\right]\exp{\left(- \frac{\lvert x_i-z\rvert^2}{2T}\right)}\\
                &\geq\exp{\left(- \frac{\lvert x_i-z\rvert^2}{2T}\right)}\mathbb{P}_0(A).
        \end{aligned}
    \end{equation}
    Let $\sigma=\inf \{t\geq0: \lvert B(t)\rvert \geq  1\}$. Using the scaling property, it holds that
    \begin{equation}
        \mathbb{P}_0(A) = \mathbb{P}_0(\sigma>T/\epsilon^2)>0.
    \end{equation}
    Hence we may choose $q_-=e^{-(L^2 + 9/4)/2}\,\mathbb{P}_0(A)$ from which the claim follows.    
\end{proof}
The next Lemma achieves the second and third step of the outline: Observe that the Brownian motion alternates between visits to $Z_1$ and excursions from $\partial Z_1$ back to $\partial Z_0$. A resonator can only be hit during the time it spends inside $Z_1$, that is during a visit. This allows us to derive tail estimates for the durations of visits and excursions.
\begin{lemma}\label{lem: Visits}
    Let $z\in\Omega$ with $\lvert z_\perp\rvert\leq2$. Let $\beta_0=0$, $\alpha_1=\inf\{t\geq0:\lvert B_\perp(t)\rvert\geq2\}$ (the time $\alpha$ of \cref{lem: Killing}) and, for $k\geq1$,
    \[
        \beta_k=\inf\{t\geq\alpha_k:\lvert B_\perp(t)\rvert\leq1\},\qquad\alpha_{k+1}=\inf\{t\geq\beta_k:\lvert B_\perp(t)\rvert\geq2\},
    \]
    and write $(a)_+=\max\{a,0\}$. Let $U_k=\alpha_k-\beta_{k-1}$ (the $k$-th visit to $Z_1$) and $Y_k=\beta_k-\alpha_k$ (the $k$-th excursion out of $Z_0$). Then:
    \begin{enumerate}[label=(\alph*)]
        \item The $\alpha_k,\beta_k$ are stopping times; $\mathbb{P}_z$-almost surely they are finite, $0=\beta_0\leq\alpha_1<\beta_1<\alpha_2<\dots$, and $\alpha_k\to\infty$.
        \item $\tau\notin[\alpha_k,\beta_k]$ for all $k\geq1$.
        \item $\mathbb{P}_z(\tau>\alpha_k)\leq(1-q_-)^{k-1}$ for all $k\geq1$; in particular $\tau<\infty$ $\mathbb{P}_z$-almost surely.
        \item For all $k\geq1$ and $s>0$,
        \begin{align}
            \mathbb{P}_z\left(\tau>\beta_{k-1},\ U_k>s\right)&\leq(1-q_-)^{(k-2)_+}\min\{1,36/s\},\label{eq: VisitDuration}\\
            \mathbb{P}_z\left(\tau>\alpha_k,\ Y_k>s\right)&\leq(1-q_-)^{k-1}\,12/\log(2+s).\label{eq: ExcursionDuration}
        \end{align}
    \end{enumerate}
\end{lemma}
\begin{proof}
    (a) The $\alpha_k,\beta_k$ are successive hitting times of the closed sets $\mathbb{R}^3\setminus Z_1$ and $\overline{Z_0}$, hence stopping times \cite[Proposition 1.15]{chung2012brownian}. Moreover they are functions of $B_\perp$ only. Whenever finite, $\lvert B_\perp(\alpha_k)\rvert=2$ and $\lvert B_\perp(\beta_k)\rvert=1$ by continuity, which gives the strict inequalities. Finiteness follows inductively from \cref{lem: DiskExit,lem: Excursion} and \eqref{eq: StrongMarkov} at $\alpha_k$ and $\beta_k$. If $\alpha_k$ increased to a finite limit, $B_\perp$ would not be uniformly continuous on a compact interval, since $\lvert B_\perp(\alpha_k)\rvert-\lvert B_\perp(\beta_k)\rvert=1$ and $\alpha_k<\beta_k<\alpha_{k+1}$.

    (b) On $[\alpha_k,\beta_k]$ we have $\lvert B_\perp\rvert\geq1>\rho\geq\lvert a_i\rvert+r_i$ for all $i$ by \eqref{eq: Cylinder Assumption} and $\rho\leq1/2$, so $B\in\Omega$.

    (c) For $k=1$ the claim is trivial. By (b), $\{\tau>\alpha_k\}=\{\beta_k<\tau\}$ up to a null set; on this event, $\tau=\beta_k+\tau\circ\theta_{\beta_k}$ and $\alpha_{k+1}=\beta_k+\alpha_1\circ\theta_{\beta_k}$. By \eqref{eq: StrongMarkov}, $\mathbb{P}_z(\tau>\alpha_{k+1})=\mathbb{E}_z[\mathrm{1}_{\{\tau>\beta_k\}}\mathbb{P}_{B(\beta_k)}(\tau>\alpha_1)]$. For $w\in\partial Z_0$ we have $\mathbb{P}_w(\tau=\alpha_1)=0$ (as $B(\alpha_1)\in\partial Z_1\subset\Omega$), so $\mathbb{P}_w(\tau>\alpha_1)\leq1-q_-$ by \cref{lem: Killing}. As $B(\beta_k)\in\partial Z_0$, $\mathbb{P}_z(\tau>\alpha_{k+1})\leq(1-q_-)\mathbb{P}_z(\tau>\alpha_k)$, and (c) follows by induction; finally $\mathbb{P}_z(\tau=\infty)\leq\mathbb{P}_z(\tau>\alpha_k)\to0$.

    (d) $U_k$ is the exit time from the disc of radius $2$ of the planar motion restarted at $B_\perp(\beta_{k-1})$, $\lvert B_\perp(\beta_{k-1})\rvert\leq2$. By \eqref{eq: StrongMarkov} and \cref{lem: DiskExit} ($R=2$), $\mathbb{P}_z(\tau>\beta_{k-1},U_k>s)\leq\mathbb{P}_z(\tau>\beta_{k-1})\min\{1,36/s\}$; here $\mathbb{P}_z(\tau>\beta_0)\leq1$, and $\mathbb{P}_z(\tau>\beta_{k-1})=\mathbb{P}_z(\tau>\alpha_{k-1})\leq(1-q_-)^{k-2}$ for $k\geq2$ by (b), (c). Likewise, $Y_k$ is the hitting time of the closed unit disc by the planar motion restarted at $B_\perp(\alpha_k)$, $\lvert B_\perp(\alpha_k)\rvert=2$, and \eqref{eq: StrongMarkov}, \cref{lem: Excursion} and (c) give \eqref{eq: ExcursionDuration}.
\end{proof}

Using the geometric tail of the number of visits together with the duration estimates, we obtain tail estimates on the survival duration. This is step three of the outline.

\begin{proposition}[Logarithmic survival bound]\label{prop: Survival}
    Let $K_0=635\,q_-^{-3}$, which depends only on $r_{\min}$ and $L$. For all $z\in\Omega$ with $\lvert z_\perp\rvert\leq2$ and $t\geq0$,
    \begin{equation}\label{eq: Survival}
        \mathbb{P}_z(\tau>t)\leq\frac{K_0}{\log(2+t)}\qquad\text{and}\qquad\mathbb{E}_z[\tau\wedge t]\leq\frac{3K_0\,t}{\log(2+t)} .
    \end{equation}
\end{proposition}
\begin{proof}
    We begin with the tail estimate on $\mathbb{P}_z(\tau >t)$. Let $w_k=\frac{1}{k(k+1)}$, so that $\sum_{k\geq1}w_k=1$. Almost surely $\tau<\infty$, and since $\alpha_k\to\infty$ we may define $G=\min\{k\geq1:\tau<\alpha_k\}$. As $\tau<\alpha_G=\sum_{k=1}^GU_k+\sum_{k=1}^{G-1}Y_k$, on $\{\tau>t\}$ it is impossible that $U_k\leq w_kt/2$ for all $k\leq G$ and $Y_k\leq w_kt/2$ for all $k\leq G-1$. Hence, as $\{k\leq G\}\subset\{\tau>\beta_{k-1}\}$ and $\{k\leq G-1\}\subset\{\tau>\alpha_k\}$  \cref{lem: Visits}(b),
    by \cref{lem: Visits}(d),
    \begin{equation}\label{eq: Survival1}
        \begin{aligned}
            \mathbb{P}_z(\tau>t)&\leq \mathbb{P}_z(\alpha_G>t)\\
            &\leq\mathbb{P}_z\left(\bigcup_{k\geq 1}\{\tau>\alpha_k,Y_k>w_k t/2\}\cup\{\tau>\beta_{k-1},U_k>w_k t/2\}\right)\\
            &\leq \sum_{k\geq1}(1-q_-)^{(k-2)_+}\min\left\{1,\frac{72k(k+1)}{t}\right\}\\
        &+\sum_{k\geq1}(1-q_-)^{k-1}\min\left\{1,\frac{12}{\log\left(2+\frac{t}{2k(k+1)}\right)}\right\}.
        \end{aligned}
    \end{equation}
    We thus have to upper bound these two sums. Let $t\geq512$. Using 
    \begin{equation}
        \sum_{j\geq0}(j+1)(j+2)x^j=2(1-x)^{-3}\qquad\text{and}\qquad(j+2)(j+3)\leq3(j+1)(j+2),
    \end{equation}
     we get
     \begin{equation}
         \sum_{k\geq1}k(k+1)(1-q_-)^{(k-2)_+}\leq2+6q_-^{-3}\leq8q_-^{-3}.
     \end{equation}     
     So the first sum in \cref{eq: Survival1}
     is at most $576/(q_-^3t)\leq576/(q_-^3\log(2+t))$.
     
     In the second sum let $k_t=\lfloor t^{1/4}\rfloor$. For $k\leq k_t$ we have $2k(k+1)\leq4\sqrt t$ and hence $\log(2+\frac{t}{2k(k+1)})\geq\frac12\log t-\log4\geq\frac14\log(2t)\geq\frac14\log(2+t)$, these terms contribute at most $48/(q_-\log(2+t))$. The terms with $k>k_t$ contribute at most
    \[
        \sum_{k>k_t}(1-q_-)^{k-1}\leq\frac{e^{-q_-(t^{1/4}-1)}}{q_-}\leq\frac{4e}{q_-^2\log(2+t)},
    \]
    where we used $\log(2+t)\leq4\log(2+t^{1/4})\leq4(1+t^{1/4})$ and $(1+y)e^{-q_-y}\leq q_-^{-1}$ for $y\geq0$. Altogether we therefore obtain
    \begin{equation}
        \begin{aligned}
            \mathbb{P}_z(\tau>t)&\leq \frac{576}{q_-^3\log(2+t)} + \frac{48}{q_-\log(2+t)} + \frac{4e}{q_-^2\log(2+t)}\\
            &\leq \frac{635q_-^{-3}}{\log(2+t)},
        \end{aligned}
    \end{equation}
    which proves the first claim for $t\geq512$. For $t<512$ it holds trivially, as $635/\log(514)>1$.
  
    The second bound now follows by integrating over the tail estimate on $\mathbb{P}_z(\tau>t)$. That is, for $t>4$ we have
    \begin{equation}        
        \begin{aligned}
            \mathbb{E}_z(\tau\wedge t)&=\int_0^t\mathbb{P}_z(\tau>s)ds\\
            &\leq \sqrt{t} + \int_{\sqrt{t}}^t\frac{K_0}{\log(2+s)}ds\\
            &\leq \sqrt{t} + \frac{K_0t}{\log(2+\sqrt{t})}\\
            &\leq\frac{3K_0t}{\log(2+t)},
        \end{aligned}
    \end{equation}
    where in the last step we used $(2+\sqrt t)^2\geq2+t$ and $\sqrt{t}\leq t/\log(2+t)\leq K_0t/\log(2+t)$, as $\log(2+t)\leq\sqrt{t}$ for $t>4$. For $t\leq4$ the bound follows from $\mathbb{E}_z(\tau\wedge t)\leq t$ and $\log(2+t)\leq\log6<3K_0$.
\end{proof}
\subsubsection{Proof of the upper bound in \cref{thm: ChainDecay}}

By \cref{lem: Flux} it is sufficient to upper bound $V_i$ on a shell around $D_j$. \cref{lem:occupation} turns this into an estimate of the occupation time $\mathcal{G}^\Omega(z,A_{i}^\delta)$. To this end, fix $\delta=\min\{d_{\min}/2,1\}$, let $z\in\partial B_{r_j+\delta}(x_j)$ and $y\in A_i^\delta$. Since $\rho\leq1/2$, we have $\lvert z_\perp\rvert,\lvert y_\perp\rvert\leq\rho+\delta\leq2$, and \eqref{eq: MinimalDistance} yields $\lvert z-x_i\rvert\geq d_{ij}-r_j-\delta\geq r_i+\delta$, so that \cref{lem:occupation} and \cref{prop: Survival} apply. Rewriting as in \cref{eq: OccupationTime} allows us to make use of \cref{lem:split}. Concretely, we apply it with $\Pi=\Pi_\parallel$: let $x,y\in\Omega$ with $\lvert x_\perp\rvert,\lvert y_\perp\rvert\leq2$ and $R=\lvert\Pi_\parallel(x-y)\rvert>0$. We then have
\begin{equation}
    \begin{aligned}
        p_t^\Omega(x,y)&\le (2/\pi {t})^{3/2}
        [\mathbb{P}_y(\tau>t/4)\mathbb{P}_x\left(\tau>t/2,\ |\Pi_\parallel(B(t/2)-x)|\ge R/2\right)\\
        &+\mathbb{P}_x(\tau>t/4)\mathbb{P}_y\left(\tau>t/2,\ |\Pi_\parallel(B(t/2)-y)|\ge R/2\right)].
    \end{aligned}
\end{equation}
To bound $\mathbb{P}_x\left(\tau>t/2,\ |\Pi_\parallel(B(t/2)-x)|\ge R/2\right)$, let $T=\tau\wedge t/2$ and set $\mathcal{M}(t)=\lvert \Pi_\parallel(B(t)-x)\rvert^2-t$. This is a martingale; applying the optional stopping theorem \cite[Proposition 2.42]{morters2010brownian} to $T\leq t/2$ (admissible, as $\sup_{u\leq t/2}\lvert\mathcal{M}(u)\rvert$ is integrable by Doob's inequality), we therefore have
\begin{equation}
    \mathbb{E}_x\left[\mathcal{M}(T)\right]=\mathbb{E}_x\left[\mathcal{M}(0)\right]=0.
\end{equation}
A Markov inequality with \cref{prop: Survival} thus yields
\begin{equation}
    \mathbb{P}_x\left(\tau>t/2,\ |\Pi_\parallel(B(t/2)-x)|\ge R/2\right)\leq \frac{6K_0t}{R^2\log(2 + t/2)}.
\end{equation}
The above bound together with the last part of \cref{lem:split} thus yields
\begin{equation}
    \begin{aligned}
        p_t^\Omega(x,y)&\leq 12(2 / \pi t)^{3/2}\frac{K_0}{\log(2 + t/4)}\frac{K_0 t}{R^2\log(2 + t/2)}\\
        &\leq12(2 / \pi t)^{3/2}\frac{K_0^2 t}{R^2\log^{2}(2 + t/4)}.
    \end{aligned}
\end{equation}
Alternatively, we may apply \cref{prop: Survival} together with the semi-group property from \cref{thm: KilledDensity} to obtain
\begin{equation}
    p_t^\Omega(x,y)\leq(2 / \pi )^{3/2}\frac{K_0^2 }{t^{3/2}\log^{2}(2 + t/4)}\min\left(1,\frac{12t}{R^2}\right).
\end{equation}
Let now $t^*=R^2/\log^4 R$ and assume $R\geq e^6$, we then have
\begin{equation}
    \begin{aligned}
        \mathcal{G}^\Omega(x,y)
        &=\int_0^\infty p_t^\Omega(x,y)dt\\
        &=\int_0^{t^*} p_t^\Omega(x,y)dt + \int_{t^*}^{R^2} p_t^\Omega(x,y)dt+\int_{R^2}^\infty p_t^\Omega(x,y)dt.
    \end{aligned}
\end{equation}
For the first integral, we have
\begin{equation}
    \begin{aligned}
        \int_0^{t^*} p_t^\Omega(x,y)dt&\leq\int_0^{t^*}(2 / \pi )^{3/2}\frac{K_0^2 }{\sqrt{t}\log^2(2 + t/4)}\frac{12}{R^2}dt\\
        &\leq(2 / \pi )^{3/2}\frac{24\sqrt{t^*}K_0^2}{R^2\log^22}\\
        &=(2 / \pi )^{3/2}\frac{24K_0^2}{R(\log R\log2)^2}.
    \end{aligned}
\end{equation}
The second integral yields
\begin{equation}
    \begin{aligned}
        \int_{t^*}^{R^2} p_t^\Omega(x,y)dt&\leq (2 / \pi )^{3/2}\frac{12K_0^2}{R^2\log^2(2 + t^*/4)}\int_{t^*}^{R^2}\frac{dt}{\sqrt{t}}
        \leq(2 / \pi )^{3/2}\frac{24K_0^2}{R\log^2(2 + t^*/4)}\\
        &\leq(2 / \pi )^{3/2}\frac{96K_0^2}{R\log^2 R},
    \end{aligned}
\end{equation}
where we used that $\log(2+t^*/4)\geq2\log R-4\log\log R-\log 4\geq\frac{1}{2}\log R$ for $R\geq e^6$.
For the third one, we obtain
\begin{equation}
    \begin{aligned}
        \int_{R^2}^\infty p_t^\Omega(x,y)dt&\leq(2 / \pi )^{3/2}\frac{2K_0^2}{R\log ^2(2 + R^2/4)}.
    \end{aligned}
\end{equation}
Putting everything together, we have
\begin{equation}
    \begin{aligned}
        \mathcal{G}^\Omega(x,y)&\leq (2 / \pi )^{3/2}\left(\frac{24K_0^2}{R(\log R\log2)^2} + \frac{96K_0^2}{R\log^2 R} +\frac{2K_0^2}{R\log ^2(2 + R^2/4)}\right)\\
        &\leq \frac{CK_0^2}{R\log^2 R},
    \end{aligned}
\end{equation}
with an absolute constant $C$. Integrating over $A_{i}^\delta$, we therefore have
\begin{equation}
    V_i(z)\leq\frac{12}{\delta^2}\mathcal{G}^\Omega(z,A_{i}^\delta)\leq\frac{C'K_0^2}{\delta^2d_{ij}\log^2(2+d_{ij})},\qquad d_{ij}\geq 2e^6,
\end{equation}
where we used that $R=\lvert\Pi_\parallel(z-y)\rvert\geq\lvert s_i-s_j\rvert-(r_i+\delta)-(r_j+\delta)\geq d_{ij}-4\geq d_{ij}/2\geq e^6$ for every $y\in A_i^\delta$, as $\lvert a_i-a_j\rvert\leq2\rho\leq1$ and $r_k+\delta\leq3/2$, that $R\mapsto R\log^2R$ is increasing, and that $\log(d_{ij}/2)\geq\frac{1}{2}\log(2+d_{ij})$ for $d_{ij}\geq8$. Taking the supremum over $z\in\partial B_{r_j+\delta}(x_j)$ in \eqref{eq: CapacitanceByPotential} yields the upper bound for $d_{ij}\geq2e^6$. For $d_{ij}<2e^6$ it follows from $V_i\leq1$ and \eqref{eq: CapacitanceByPotential}, as $d\mapsto d\log^2(2+d)$ is bounded on bounded sets.
\subsection{Screens}
The upper bounds for screens follows the same outline of the decay estimate for chains. We begin by defining the slabs
\begin{equation}
    S_k:=\{z\in\R^3: \lvert{z_\parallel}\rvert<2^k\}.
\end{equation}
As before, a uniform killing bound holds:
\begin{lemma}\label{lem: KillingScreen}
    Let $\alpha=\inf\{t\geq0:\lvert {B_\parallel}(t)\rvert\geq2\}$ be the exit time from $S_1$. For every $z\in\partial S_0$,
    \begin{equation}\label{eq: qminusScreen}
        \mathbb{P}_z(\tau<\alpha)\geq q_-\in(0,1).
    \end{equation}
    The constant $q_-$ depends only on $r_{\min}$ and $L$.
\end{lemma}
\begin{proof}
    The proof follows mutatis mutandis from \cref{lem: Killing}.
\end{proof}
As the reader may expect, the next step is to adapt \cref{lem: Visits}. 
\begin{lemma}\label{lem: VisitsScreen}
    Let $z\in\Omega$ with $\lvert {z_\parallel}\rvert\leq2$. Let $\beta_0=0$, $\alpha_1=\inf\{t\geq0:\lvert {B_\parallel}(t)\rvert\geq2\}$ (the time $\alpha$ of \cref{lem: KillingScreen}) and, for $k\geq1$,
    \[
        \beta_k=\inf\{t\geq\alpha_k:\lvert {B_\parallel}(t)\rvert\leq1\},\qquad\alpha_{k+1}=\inf\{t\geq\beta_k:\lvert {B_\parallel}(t)\rvert\geq2\},
    \]
    and write $(a)_+=\max\{a,0\}$. Let $U_k=\alpha_k-\beta_{k-1}$ (the $k$-th visit to $S_1$) and $Y_k=\beta_k-\alpha_k$ (the $k$-th excursion out of $S_0$). Then:
    \begin{enumerate}[label=(\alph*)]
        \item The $\alpha_k,\beta_k$ are stopping times; $\mathbb{P}_z$-almost surely they are finite, $0=\beta_0\leq\alpha_1<\beta_1<\alpha_2<\dots$, and $\alpha_k\to\infty$.
        \item $\tau\notin[\alpha_k,\beta_k]$ for all $k\geq1$.
        \item $\mathbb{P}_z(\tau>\alpha_k)\leq(1-q_-)^{k-1}$ for all $k\geq1$; in particular $\tau<\infty$ $\mathbb{P}_z$-almost surely.
        \item For all $k\geq1$ and $s>0$,
        \begin{align}
            \mathbb{P}_z\left(\tau>\beta_{k-1},\ U_k>s\right)&\leq(1-q_-)^{(k-2)_+}\min\{1,4/s\},\label{eq: VisitDurationScreen}\\
            \mathbb{P}_z\left(\tau>\alpha_k,\ Y_k>s\right)&\leq(1-q_-)^{k-1}\,\min\{1,s^{-1/2}\}.\label{eq: ExcursionDurationScreen}
        \end{align}
    \end{enumerate}
\end{lemma}
\begin{proof}
    Items (a)-(c) follow from the same arguments as in \cref{lem: Visits}, where in (b) we now use $\lvert B_\parallel\rvert\geq1>\rho\geq\lvert s_i\rvert+r_i$ on $[\alpha_k,\beta_k]$, by \eqref{eq: Slab Assumption} and $\rho\leq1/2$. It thus remains to prove (d). We now note that $U_k$ is the exit time from $(-2,2)$ of the one dimensional Brownian motion $B_\parallel$ started at $B_\parallel(\beta_{k-1})\in[-2,2]$. Let $u<0<v$, if we denote by $T$ the exit time of standard scalar Brownian motion from the interval $[u,v]$, it holds that $\mathbb{E}_0[T]=\lvert {u}\rvert {v}$ \cite[Theorem 2.49]{morters2010brownian}, so that the expected exit time from $(-2,2)$ of a Brownian motion started at $c\in[-2,2]$ equals $(2-c)(2+c)\leq4$. The strong Markov property together with an application of the Markov inequality thus yields
    \begin{align}
        \mathbb{P}_{z}(\tau>\beta_{k-1},U_k>s)\leq(1-q_-)^{(k-2)_+}\min(1,4/s).
    \end{align}
    Similarly, $Y_k$ can be understood as the hitting time of a point at distance $1$ by a linear Brownian motion $W$ started at $0$. By the reflection principle \cite[Theorem 2.21]{morters2010brownian}, $\mathbb{P}(\max_{[0,s]}W<1)=\mathbb{P}(\lvert W(s)\rvert<1)\leq2/\sqrt{2\pi s}\leq s^{-1/2}$, hence
    \begin{equation}
        \mathbb{P}_{z}(\tau>\alpha_{k},Y_k>s)\leq(1-q_-)^{k-1}\min\left(1,s^{-1/2}\right).
    \end{equation}
\end{proof}

\begin{proposition}[Logarithmic survival bound]\label{prop: SurvivalScreen}
    For all $z\in\Omega$ with $\lvert {z_\parallel}\rvert\leq2$ and $t\geq0$,
    \begin{equation}\label{eq: SurvivalScreen}
        \mathbb{P}_z(\tau>t)\leq 68q_-^{-3}\min\left(1,t^{-1/2}\right)\qquad\text{and}\qquad\mathbb{E}_z[(\tau\wedge t)^2]\leq\frac{272t^{3/2}}{3q_-^3} .
    \end{equation}
\end{proposition}
\begin{proof}
    The proof is similar to that of \cref{prop: Survival}. Letting again $w_k = \frac{1}{k(k+1)}$ and using \cref{lem: VisitsScreen}, we obtain
    \begin{equation}\label{eq: SurvivalScreen1}
        \begin{aligned}
            \mathbb{P}_z(\tau>t)&\leq \mathbb{P}_z(\alpha_G>t)\\
            &\leq\mathbb{P}_z\left(\bigcup_{k\geq 1}\{\tau>\alpha_k,Y_k>w_k t/2\}\cup\{\tau>\beta_{k-1},U_k>w_k t/2\}\right)\\
            &\leq \sum_{k\geq1}(1-q_-)^{(k-2)_+}\min\left\{1,\frac{8k(k+1)}{t}\right\}\\
        &+\sum_{k\geq1}(1-q_-)^{k-1}\min\left\{1,\left(\frac{2k(k+1)}{t}\right)^{1/2}\right\}.
        \end{aligned}
    \end{equation}
    An argument as in \cref{prop: Survival} yields
    \begin{equation}
        \sum_{k\geq1}(1-q_-)^{(k-2)_+}\min\left\{1,\frac{8k(k+1)}{t}\right\}\leq64q_-^{-3}\min\left(1,t^{-1}\right),
    \end{equation}
    while the second sum can be bounded by
    \begin{equation}
        \sum_{k\geq1}(1-q_-)^{k-1}\min\left\{1,\left(\frac{2k(k+1)}{t}\right)^{1/2}\right\}\leq\frac{4}{q_-^2}\min\left(1,t^{-1/2}\right).
    \end{equation}
    Putting these two together, we have
    \begin{equation}
        \mathbb{P}_z(\tau>t)\leq 68q_-^{-3}\min\left(1,t^{-1/2}\right).
    \end{equation}
    Integrating this tail bound, we obtain
    \begin{equation}
        \begin{aligned}
            \mathbb{E}_z[(\tau\wedge t)^2]&=2\int_0^ts\mathbb{P}_z(\tau>s)ds\\
            &\leq136q_-^{-3}\int_0^ts^{1/2}ds\\
            &\leq\frac{272t^{3/2}}{3q_-^3}.
        \end{aligned}
    \end{equation}
\end{proof}
\subsubsection{Proof of the upper bound in \cref{thm: ScreenDecay}}
The strategy remains the same as in the upper bound of \cref{thm: ChainDecay}. Let again $\delta=\min\{d_{\min}/2,1\}$, and let $x,y\in\Omega$ with $\lvert x_\parallel\rvert,\lvert y_\parallel\rvert\leq2$ and $R=\lvert\Pi_\perp(x-y)\rvert>0$. By \cref{lem:split} (applied with $\Pi=\Pi_\perp$ and $p=4$) we thus have
\begin{equation}\label{eq: ScreenDecay1}
    \begin{aligned}
        p_t^\Omega(x,y)&\le (2 / \pi t)^{3/2}\frac{16}{R^4}
        [\mathbb{P}_y(\tau>t/4)\mathbb{E}_x\left[\lvert\Pi_\perp(B(\tau\wedge t/2)-x)\rvert^4\right]\\
        &+\mathbb{P}_x(\tau>t/4)\mathbb{E}_y\left[\lvert\Pi_\perp(B(\tau\wedge t/2)-y)\rvert^4\right]].
    \end{aligned}
\end{equation}
Once again, noting that $W(t)=\Pi_\perp(B(t)-x)$ is a planar Brownian motion under $\mathbb{P}_x$ and that
\begin{equation}
    \mathcal{M}(t)=\lvert{W(t)}\rvert^4 - 8t\lvert W(t)\rvert^2+8t^2
\end{equation}
is a martingale, applying the optional stopping Theorem \cite[Proposition 2.42]{morters2010brownian} to the bounded stopping time $T=\tau\wedge t/2$ (admissible, as $\sup_{u\leq t/2}\lvert\mathcal{M}(u)\rvert$ is integrable by Doob's inequality) yields $\mathbb{E}_x[\mathcal{M}(T)]=0$, hence 
\begin{equation}
    \mathbb{E}_x\left[\lvert W(T)\rvert^4\right]=8\,\mathbb{E}_x\left[T\lvert W(T)\rvert^2\right]-8\,\mathbb{E}_x\left[T^2\right]\leq8\,\mathbb{E}_x\left[T^2\right]^{1/2}\mathbb{E}_x\left[\lvert W(T)\rvert^4\right]^{1/2},
\end{equation}
that is, $\mathbb{E}_x\left[\lvert\Pi_\perp(B(\tau\wedge t/2)-x)\rvert^4\right]\leq64\,\mathbb{E}_x[(\tau\wedge t/2)^2]$. Combining this, \cref{prop: SurvivalScreen} and Equation \eqref{eq: ScreenDecay1}, we obtain
\begin{equation}\label{eq: ScreenDecay2}
    p_t^\Omega(x,y)\leq(2/\pi t)^{3/2}\frac{16}{R^4}\cdot2\cdot68q_-^{-3}\min\left(1,2t^{-1/2}\right)\cdot64\cdot\frac{272}{3q_-^3}\left(\frac{t}{2}\right)^{3/2}\leq\frac{K_1}{R^4}\min\left(1,t^{-1/2}\right),
\end{equation}
where $K_1=4.6\cdot10^6q_-^{-6}$. 
On the other hand, the semigroup property from \cref{thm: KilledDensity} and \cref{prop: SurvivalScreen} give, as in the proof of the upper bound of \cref{thm: ChainDecay},
\begin{equation}\label{eq: ScreenDecay3}
    p_t^\Omega(x,y)\leq(2/\pi t)^{3/2}\mathbb{P}_y(\tau>t/4)\mathbb{P}_x(\tau>t/2)\leq\frac{K_2}{t^{5/2}},\qquad t\geq4,
\end{equation}
with $K_2=6644q_-^{-6}$. Let now $R\geq2$. Using \eqref{eq: ScreenDecay2} on $[0,R^2]$ and \eqref{eq: ScreenDecay3} on $[R^2,\infty)$, we obtain
\begin{equation}
    \mathcal{G}^\Omega(x,y)=\int_0^\infty p_t^\Omega(x,y)dt\leq\frac{K_1}{R^4}\left(1+\int_1^{R^2}\frac{dt}{\sqrt{t}}\right)+K_2\int_{R^2}^\infty\frac{dt}{t^{5/2}}\leq\frac{2K_1}{R^3}+\frac{2K_2}{3R^3}\leq\frac{3K_1}{R^3}.
\end{equation}
Let finally $z\in\partial B_{r_j+\delta}(x_j)$ and $y\in A_i^\delta$. Then $\lvert z_\parallel\rvert,\lvert y_\parallel\rvert\leq\rho+\delta\leq2$, \cref{lem:occupation} applies as in the proof of \cref{thm: ChainDecay}, and $R=\lvert\Pi_\perp(z-y)\rvert\geq\lvert a_i-a_j\rvert-(r_i+\delta)-(r_j+\delta)\geq d_{ij}-4\geq d_{ij}/2$ for $d_{ij}\geq8$, as $\lvert s_i-s_j\rvert\leq2\rho\leq1$ and $r_k+\delta\leq3/2$.
Finally, integrating over $A_{i}^\delta$, we obtain
\begin{equation}
    V_i(z)\leq \frac{12\lvert A_{i}^\delta\rvert\, 3K_1}{\delta^2R^3}\leq\frac{288\lvert A_i^\delta\rvert K_1}{\delta^2d_{ij}^3},\qquad d_{ij}\geq8,
\end{equation}
which proves the claim for $d_{ij}\geq8$ after taking the supremum over $z$ in \eqref{eq: CapacitanceByPotential}; for $d_{ij}<8$ it follows from $V_i\leq1$ and \eqref{eq: CapacitanceByPotential}.
\subsection{Crystals}
Proving the upper bound for crystals follows a slightly different approach. The assumptions of \cref{thm: CrystalDecay} imply that at any point $z \in \Omega$, there exists some $q_-$ such that 
\begin{equation}
    \mathbb{P}_z(\tau <1)\geq q_-\in(0,1).
\end{equation}
This differs from \cref{lem: Killing,lem: KillingScreen} in that we do not require $z$ to lie inside some slab or cylinder. Consequently, this estimate will be enough to force exponential decay.
\begin{lemma}\label{lem: KillingCrystal}
    Under the assumptions of \cref{thm: CrystalDecay}, there exists some $q_->0$ such that for every $z \in \Omega$
    \begin{equation}\label{eq: qminusCrystal}
        \mathbb{P}_z(\tau<1)\geq q_-.
    \end{equation}
    The constant $q_-$ depends only on $r_{\min}$ and $L$.
\end{lemma}
\begin{proof}
    The proof follows mutatis mutandis from \cref{lem: Killing}.
\end{proof}

\begin{proposition}\label{prop: SurvivalCrystal}
    There exist $C,\mu>0$ such that for all $z\in\Omega$
    \begin{equation}\label{eq: SurvivalCrystal}
        \mathbb{P}_z(\tau>t)\leq Ce^{-\mu t}.
    \end{equation}
\end{proposition}
\begin{proof}
    Let $E_n=\{\tau > n\}$, the strong Markov property then yields
    \begin{equation}
        \mathbb{P}_z(E_n)=\mathbb{E}_z\left[\mathrm{1}_{E_{n-1}}\mathbb{P}_{B(n-1)}(\tau>1)\right]\leq (1-q_-)\mathbb{P}_z(E_{n-1}).
    \end{equation}
    By induction, we therefore have
    \begin{equation}
        \mathbb{P}_z(\tau>n)\leq (1-q_-)^n,
    \end{equation}
    or more generally,
    \begin{equation}
        \mathbb{P}_z(\tau>t)\leq (1-q_-)^{t-1}.
    \end{equation}
    Setting $\mu = -\log(1-q_-)$ and $C=(1-q_-)^{-1}$, the claim follows.
\end{proof}
\subsubsection{Proof of the upper bound in \cref{thm: CrystalDecay}}
Once again we want to obtain an upper bound on the killed transition density $p_t^\Omega$. Thanks to the exponential decay bound of \cref{prop: SurvivalCrystal}, the estimate is much less involved: The semigroup property of \cref{thm: KilledDensity} together with $p^\Omega_{t/2}\leq p_{t/2}\leq(\pi t)^{-3/2}$ and \cref{prop: SurvivalCrystal} yield
\begin{equation}
    p_t^\Omega(x,y)\leq C\pi^{-3/2}t^{-3/2}e^{-\mu t/2}.
\end{equation}
On the other hand, we always have $p_t^\Omega \leq p_t$. Taking the geometric mean of these two upper bounds, we thus have
\begin{equation}
    p_t^\Omega(x,y)\leq C t^{-3/2}\exp(-\lvert x-y\rvert^2/4t - \mu t/4).
\end{equation}
Define
\begin{equation}
    J(c) = \int_0^\infty \exp(-u^2-c^2/u^2)du,
\end{equation}
by substituting $u=\lvert x-y\rvert/(2\sqrt{t})$ one can show that
\begin{equation}
   \mathcal{G}^\Omega(x,y)= \int_0^\infty p_t^\Omega(x,y)dt\leq \frac{4C}{\lvert x - y\rvert}J\left(\frac{\lvert x-y\rvert\sqrt{\mu}}{4}\right).
\end{equation}
Differentiating under the integral sign and substituting $v=c/u$ gives
\begin{equation}
    J'(c) = -2J(c).
\end{equation}
Simultaneously, we have $J(0)=\sqrt{\pi}/2$ and hence $ J(c) = \frac{\sqrt{\pi}}{2}e^{-2c}$. We therefore have
\begin{equation}
   \mathcal{G}^\Omega(x,y) \leq \frac{C}{\lvert x-y\rvert}\exp(-\lvert x - y\rvert\sqrt{\mu}/2).
\end{equation}
Integrating over $A_{i}^\delta$ with $\delta=\min\{d_{\min}/2,1\}$, where $\lvert z-y\rvert\geq d_{ij}-2r_{\max}-2$ for $z\in\partial B_{r_j+\delta}(x_j)$ and $y\in A_i^\delta$, together with \cref{lem:occupation} and \eqref{eq: CapacitanceByPotential} yields the claim.

\section{Lower bounds}\label{sec: Lower}

Throughout this section, we may assume that $d_{\min}\leq1$,
replacing the lower separation bound in
\eqref{eq: MinimalDistance} by $\min\{d_{\min},1\}$ if necessary.
This does not restrict the admissible configurations, and the
constants obtained below still depend only on the geometric
parameters in the statements of the main theorems.
For chains and screens, we also retain the normalisation
$\rho\leq1/2$ provided by \cref{lem: Scaling}.

Finally, we will derive the lower bounds of \cref{thm: ChainDecay,thm: ScreenDecay} and \cref{thm: CrystalDecay}. Once again, the arguments in the chain and screen setting follow the same strategy:
\begin{enumerate}
    \item We show that there exists a small ball $Q$ in $Z_1^c$, respectively $S_1^c$, on which $V_i$ is uniformly lower bounded.
    \item Using this lower bound on $Q$, we derive a lower bound for all $x$ in some subset $U\subset\Omega$ in terms of $\mathbb{P}_x(\tau_Q<\tau_U)$.
    \item Using \cref{lem:hit-lower}, we lower bound $\mathbb{P}_x(\tau_Q<\tau_U)$ in terms of $\mathcal{G}^U\mathrm{1}_Q(x)$.
    \item The respective geometries then yield a lower bound on $\mathcal{G}^U$. 
    \item Finally, we turn the lower bound on $V_i$ into a lower bound for $\lvert \mathcal{C}_{ij}\rvert$.
\end{enumerate}
We will develop this machinery in detail in Subsection \ref{ssec: ChainLower}, afterwards, in Subsection \ref{ssec: ScreenLower} we will be a bit laissez-faire with the arguments as they are so similar. For crystals, a shorter argument can be made using the Cameron-Martin Theorem.
\subsection{Chains}\label{ssec: ChainLower}
As a first step, we show that to each resonator $D_i$ there exists some ball $Q_i\subset\Omega$, on which $V_i$ is uniformly lower bounded.
\begin{lemma}\label{lem: PotentialLower}
    Let $n\in\R^2\times\{0\}$ have unit norm. Set $z_i = x_i + 3n$ and $Q_i=\overline{B_{1/2}(z_i)}$. Then we have $Q_i\subset{Z_1^c\subset}\Omega$. Moreover, there exists a constant $c_0>0$, depending only on $r_{\min}$ and $d_{\min}$ in the normalised configuration, and independent of $i$, $n$, and the particular configuration such that
    \begin{equation}
        V_i\mid_{Q_i}\geq c_0.
    \end{equation}
\end{lemma}
\begin{proof}
Fix $q\in Q_i$. Since $|a_i|+r_i\leq1/2$,
\[
    |q_\perp|
    \geq3-|a_i|-1/2
    \geq2+r_i>2.
\]
Thus $Q_i\subset Z_1^c\subset\Omega$, and every point of
$Q_i$ has distance greater than $3/2$ from $\mathcal{D}$.

For each $k\neq i$, consider the enlarged ball
$B_{r_k+d_{\min}/3}(x_k)$. Their closures are pairwise
disjoint. Moreover, both $q$ and $x_i$ lie outside all these
balls. The segment joining $q$ to $x_i$ meets only finitely
many of them, by local finiteness. Replace each portion of
the segment inside an enlarged ball by a shortest
great-circle arc on its boundary.

The resulting rectifiable path joins $q$ to $x_i$ and
remains at distance at least $d_{\min}/3$ from every
$D_k$ with $k\neq i$. A shortest spherical arc has length
at most $\pi/2$ times the length of its chord. Since the
replaced portions of the segment are disjoint, the length
$\ell$ of the resulting path satisfies
\[
    \ell\leq\frac{\pi}{2}|q-x_i|
    \leq\frac{7\pi}{4}=:L_*.
\]
Parametrise this path at constant speed on $[0,1]$, and
denote the parametrisation by $\widetilde{\gamma}$.
Then
\[
    \widetilde{\gamma}(0)=q,\qquad
    \widetilde{\gamma}(1)=x_i,\qquad
    \int_0^1|\widetilde{\gamma}'(t)|^2\,dt
    =\ell^2\leq L_*^2.
\]

Set
\[
    \epsilon_*=\min\{r_{\min}/2,d_{\min}/6\},
    \qquad
    E_q=
    \left\{
      \sup_{0\leq t\leq1}
      |B(t)-\widetilde{\gamma}(t)|<\epsilon_*
    \right\}.
\]
On $E_q$, the Brownian path avoids every resonator other
than $D_i$ until time $1$, and $B(1)\in D_i$.
By continuity,
\[
    E_q\subset\{\tau<\infty,\ B(\tau)\in\partial D_i\}.
\]
The Cameron--Martin theorem, together with the symmetry
under $B\mapsto-B$ of the event
$\{\sup_{0\leq t\leq1}|B(t)|<\epsilon_*\}$, therefore gives,
as in the proof of \cref{lem: Killing},
\[
\begin{aligned}
    V_i(q)
    &\geq\mathbb{P}_q(E_q)\\
    &\geq
    \exp\left(
      -\frac12\int_0^1|\widetilde{\gamma}'(t)|^2\,dt
    \right)
    \mathbb{P}_0\left(
      \sup_{0\leq t\leq1}|B(t)|<\epsilon_*
    \right)\\
    &\geq
    e^{-L_*^2/2}
    \mathbb{P}_0\left(
      \sup_{0\leq t\leq1}|B(t)|<\epsilon_*
    \right)
    =:c_0>0.
\end{aligned}
\]
The last expression depends only on $r_{\min}$ and
$d_{\min}$ and is independent of $q$, $i$, $n$, and the
particular configuration. This proves the claim.

\end{proof}

\begin{lemma}\label{lem:hit-lower}
Let $Q=\overline B_{\varrho}(u)\subset U$ and $x\in U\setminus Q$.
For $\tau_Q=\inf\{t\ge0:B(t)\in Q\}$ and $\tau_U=\inf\{t\geq0:B(t)\notin U\}$,
\begin{equation}\label{eq:hit-lower}
 \mathbb{P}_x(\tau_Q<\tau_U)\ge \varrho^{-2}\int_Q\mathcal{G}^U(x,v)\,dv.
\end{equation}
If also $U\subset\Omega$ and $V_i\ge {m}>0$ on $Q$, then
$V_i(x)\ge {m}\mathbb{P}_x(\tau_Q<\tau_U)$.
\end{lemma}
\begin{proof}
There is no occupation of $Q$ before $\tau_Q$. Applying the strong Markov property \eqref{eq: StrongMarkov} at $\tau_Q$ with $E=\{\tau_Q<\tau_U\}\in\mathcal{F}^+(\tau_Q)$ together with domain monotonicity of the Green's function therefore yields
\begin{equation}
    \begin{aligned}
        \int_Q\mathcal{G}^U(x,y)dy &=\mathcal{G}^U\mathrm{1}_Q(x)\\
        &=\mathbb{E}_x\left[\int_0^{\tau_U}\mathrm{1}_Q(B(t))dt\right]\\
        &=\mathbb{E}_x\left[\mathrm{1}_{\tau_Q<\tau_U}\int_{\tau_Q}^{\tau_U}\mathrm{1}_Q(B(t))dt\right]\\
        &=\mathbb{E}_x\left[\mathrm{1}_{\tau_Q<\tau_U}\mathbb{E}_{B(\tau_Q)}\left[\int_0^{\tau_U}\mathrm{1}_Q(B(t))dt\right]\right]\\
        &\leq\mathbb{P}_x(\tau_Q<\tau_U)\sup_{w\in Q}\int_{Q}\mathcal{G}^0(w,y)dy.
    \end{aligned} 
\end{equation}
For every $w$, the layer-cake formula gives
\begin{align*}
 \int_Q{dv\over|w-v|}
 &=\int_0^\infty |Q\cap B_{1/s}(w)|\,ds\\
 &\le\int_0^{1/\varrho}{4\pi \varrho^3\over3}\,ds
       +\int_{1/\varrho}^\infty{4\pi\over3s^3}\,ds
 =2\pi \varrho^2.
\end{align*}
This proves \eqref{eq:hit-lower}. To obtain the last claim, note that
\begin{equation}
    \begin{aligned}
        V_i(x)&=\mathbb{E}_x[\mathrm{1}_{\tau<\infty}\mathrm{1}_{\partial D_i}(B(\tau))]\\
        &\geq \mathbb{E}_x[\mathrm{1}_{\tau_Q<\tau_U}\mathrm{1}_{\tau<\infty}\mathrm{1}_{\partial D_i}(B(\tau))]\\
        &=\mathbb{E}_x\left[\mathrm{1}_{\tau_Q<\tau_U}V_i(B(\tau_Q))\right]\\
        &\geq {m}\mathbb{P}_x(\tau_Q<\tau_U),
    \end{aligned}
\end{equation}
where we used that $\{\tau_Q<\tau_U\}\subset\{\tau_Q<\tau\}$ when $U\subset\Omega$, so that $\tau=\tau_Q+\tau\circ\theta_{\tau_Q}$ on this event and \eqref{eq: StrongMarkov} together with \eqref{eq: NewtonianPotential} give the middle equality, and that $V_i\geq m$ on $Q\ni B(\tau_Q)$.
\end{proof}
Setting $U=Z_1^c$, an application of \cref{lem:hit-lower} immediately yields
\begin{equation}
    V_i(x)\geq c_0\mathbb{P}_x(\tau_{Q_i}<\tau_{Z_1^c})\geq 4c_0\mathcal{G}^{Z_1^c}\mathrm{1}_{Q_i}(x).
\end{equation}
The next step is thus to derive an expression for $\mathcal{G}^{Z_1^c}$. To this end, we require the following result on the Green's function of planar exterior domains \cite[Theorem 1.2]{collet2000asymptotic}
\begin{theorem}\label{thm: PuncturedPlanar}
    Let $\Omega=\mathbb{R}^2\setminus\overline{B_r(0)}$ for some $r>0$. For any $x,y \in\Omega$ it holds that
    \begin{equation}
        \lim_{t\rightarrow\infty}t(\log t)^2p_t^\Omega(x,y)=\frac{2\log(\lvert x\rvert/r)\log(\lvert y\rvert/r)}{\pi}.
    \end{equation}
    The convergence is uniform on any compact subset of $\Omega\times\Omega$.
\end{theorem}
The case $r=1$ is \cite[Theorem 1.2]{collet2000asymptotic}; the general case follows from the scaling identity $p_t^\Omega(x,y)=r^{-2}p^{\Omega/r}_{t/r^2}(x/r,y/r)$ together with $(\log t)^2/(\log(t/r^2))^2\rightarrow1$.

Upper and lower bounds for Dirichlet heat kernels in the
exterior of compact sets were established by Grigor'yan
and Saloff-Coste \cite{grigoryan2002dirichlet}.
Uchiyama \cite{uchiyama2018transition} obtained uniform
large-time asymptotics for the transition density of
Brownian motion killed on a bounded set, allowing the
spatial variables to vary with time.
For our purposes, the convergence on compact subsets
stated in \cref{thm: PuncturedPlanar} is sufficient.

With this in hand, we now can prove the following Theorem
\begin{theorem}
    Fix $2<\zeta_-<\zeta_+<\infty$ and $\epsilon>0$. Then there exists some $R_\epsilon=R_\epsilon(\zeta_-,\zeta_+,\epsilon)$ such that for any $x,y \in Z_1^c$ with $\zeta_-<\lvert x_\perp\rvert,\lvert y_\perp\rvert<\zeta_+$ and $R=\lvert \Pi_\parallel(x-y)\rvert>R_\epsilon$ we have
    \begin{equation}
        \mathcal{G}^{Z_1^c}(x,y)\geq(1-\epsilon)\frac{\log(\lvert {x_\perp}\rvert/2)\log(\lvert {y_\perp}\rvert/2)}{2\pi R(\log R)^2}.
    \end{equation}
\end{theorem}
\begin{proof}
    Let $q_t(x,y)$ denote the killed density of $\overline{B_{2}(0)}^c$, so that \cref{thm: PuncturedPlanar} applies with $r=2$. Independence of axial and radial Brownian motion then implies that 
    \begin{equation}
        p_t^{Z_1^c}(x,y)=q_t(x_\perp,y_\perp)p^{(1)}_t(\Pi_\parallel(x-y)),
    \end{equation}
    where $p^{(1)}_t(u)=(2\pi t)^{-1/2}e^{-u^2/(2t)}$ denotes the one-dimensional Gaussian kernel. Suppose the claim fails. Then there is a sequence $x_n,y_n\in Z_1^c$ with $R_n=\lvert\Pi_\parallel(x_n-y_n)\rvert\rightarrow \infty$ and $\zeta_-<\lvert (x_n)_\perp\rvert,\lvert (y_n)_\perp\rvert<\zeta_+$ for which the asserted inequality fails. Writing $\xi_n=(x_n)_\perp$, $\xi_n'=(y_n)_\perp$ and $\ell_n=\log(\lvert \xi_n\rvert/2)\log(\lvert \xi_n'\rvert/2)$, \cref{thm: PuncturedPlanar} then yields
    \begin{equation}
        \lim_{n\rightarrow\infty} q_{R_{n}^2s}(\xi_n,\xi_n')\frac{2\pi R_{n}^2(\log R_{n})^2}{\ell_n} =s^{-1}\quad\text{for every }s>0.
    \end{equation}
    The scaling property  of  Brownian motion and Fatou's Lemma then yield
    \begin{equation}
        \begin{aligned}
            \liminf_{n\rightarrow\infty}\frac{2\pi R_n(\log R_n)^2}{\ell_n}\mathcal{G}^{Z_1^c}(x_n,y_n)\\
            &=\liminf_{n\rightarrow\infty}\int_0^\infty \frac{2\pi R_n(\log R_n)^2} {\ell_n}q_t(\xi_n,\xi_n')\frac{e^{-R_n^2/(2t)}}{\sqrt{2\pi t}}dt \\
            &\geq\int_0^\infty \liminf_{n\rightarrow\infty}\frac{2\pi R_n^2(\log R_n)^2}{\ell_n} q_{R_n^2s}(\xi_n,\xi_n')\frac{e^{-1/(2s)}}{\sqrt{2\pi s}}ds\\
            &\geq\int_0^\infty s^{-1}\frac{e^{-1/(2s)}}{\sqrt{2\pi s}}ds=1.
        \end{aligned}
    \end{equation}
    This contradicts the choice of the sequence, which proves the claim.
\end{proof}
Since $5/2\leq\lvert(z_j)_\perp\rvert\leq7/2$ and $9/4<\lvert y_\perp\rvert<15/4$ for $y\in B_{1/4}(z_i)\subset Q_i$, while $d_{ij}-5/4\leq R\leq d_{ij}+2$ for such $y$, for $d_{ij}\geq {R_\epsilon+2}$ we thus obtain from \cref{lem:hit-lower} and the previous theorem (with $\zeta_-=9/4$, $\zeta_+=4$)
\begin{equation}\label{eq: ChainLower1}
    V_i(z_j)\geq \frac{c_0c_1}{ (d_{ij} + 2)(\log (d_{ij} +2))^2},
\end{equation}
where $c_1$ is a constant independent of the geometric parameters.

\subsubsection*{Proof  of the lower bound in \cref{thm: ChainDecay}}\label{sssec: ChainLowerProof}
Fix $\epsilon=1/2$ in \eqref{eq: ChainLower1} and set $R_0=R_\epsilon+2$, which is an absolute constant. For $d_{ij}<R_0$ the segment from $z_j$ to $x_i$ has length at most $R_0+2r_{\max}+3$, and the construction in the proof of \cref{lem: PotentialLower}, applied to this segment, yields $V_i(z_j)\geq c_0'$ for some $c_0'>0$ depending only on $r_{\min},r_{\max},d_{\min}$. As $d_{ij}\geq d_{\min}$, and $d_{ij}+2\leq2d_{ij}$ for $d_{ij}\geq R_0$, we thus have in all cases
\begin{equation}\label{eq: ChainLower2}
    V_i(z_j)\geq\frac{c}{d_{ij}\log^2(2+d_{ij})},
\end{equation}
with $c>0$ depending only on $r_{\min},r_{\max},d_{\min}$.

Consequently, $B_{4r_h}(y)\subset\Omega$ for every point $y$
on the path. A chain of overlapping balls of radius $r_h$,
followed by repeated application of Harnack's inequality,
therefore gives
\begin{equation}
    V_i(z')\geq27^{-N}V_i(z_j),
\end{equation}
where $N$ depends only on the prescribed geometric bounds.

Every point of $\partial B_{r_j+d_{\min}/3}(x_j)$ has
distance at least $d_{\min}/3$ from $\mathcal{D}$.
Joining $z'$ to any other point of this sphere by a shortest
great-circle arc and applying the same argument gives
\begin{equation}
    \min_{z\in\partial B_{r_j+d_{\min}/3}(x_j)}V_i(z)
    \geq27^{-(N+N')}V_i(z_j),
\end{equation}
where $N'$ depends only on $r_{\max}$ and $d_{\min}$.

Finally, apply \eqref{eq: SphericalMean} with
$\delta=d_{\min}/2$, $u=V_i$, and
$\varrho=r_j+d_{\min}/3$.
Writing $\mu(\varrho)$ for the spherical mean of $V_i$,
we obtain
\begin{equation}\label{eq: LowerFlux}
    |\mathcal{C}_{ij}|
    =\frac{4\pi r_j\varrho}{\varrho-r_j}\,\mu(\varrho)
    \geq
    \frac{12\pi r_{\min}^2}{d_{\min}}\,
    27^{-(N+N')}V_i(z_j).
\end{equation}
Together with \eqref{eq: ChainLower2}, this proves the
lower bound.

\subsection{Screens}\label{ssec: ScreenLower}
The proof of the lower bound for screen geometries follows the arguments from the previous subsection. The proof of \cref{lem: PotentialLower} applies to screens
with $n=(0,0,1)$: the balls $Q_i$ then lie above the plane
$\{z_\parallel=2\}$ and have distance greater than $3/2$
from $\mathcal{D}$. The path construction and the
Cameron--Martin estimate are unchanged.
Lemma \ref{lem:hit-lower} also applies. It thus remains to choose an appropriate subdomain $U\subset\Omega$ and obtain lower bounds on its Green's function.

To this end, set $\mathcal{H}=\overline{S}_1^c\cap\mathbb{R}^2\times\mathbb{R}_+=\{z\in\R^3:z_\parallel>2\}$ and $U=\mathcal{H}$; in \cref{lem: PotentialLower} we take $n=(0,0,1)$, so that $Q_i=\overline{B_{1/2}(x_i+3n)}\subset\mathcal{H}$ by $\lvert s_i\rvert<1/2$. It is well known, see e.g. \cite{evans2010Partial} that the Green's function for a halfspace is given by
\begin{equation}
    \mathcal{G}^\mathcal{H}(x,y)=\mathcal{G}^0(x,y)-\mathcal{G}^0(x^*,y),
\end{equation}
where $x^*$ is the reflection of $x$ across the plane $\partial\mathcal{H}=\{z_\parallel=2\}$. Writing $h_x=x_\parallel-2$ and $h_y=y_\parallel-2$ for the distances of $x,y\in\mathcal{H}$ to $\partial\mathcal{H}$, we have $\lvert x^*-y\rvert^2-\lvert x-y\rvert^2=4h_xh_y$ and hence
\begin{equation}\label{eq: HalfspaceGreen}
    \mathcal{G}^\mathcal{H}(x,y)=\frac{2h_xh_y}{\pi\lvert x-y\rvert\,\lvert x^*-y\rvert\,(\lvert x-y\rvert+\lvert x^*-y\rvert)}.
\end{equation}
The following Lemma follows from a simple calculation.
\begin{lemma}\label{lem: HalfspaceGreen}
    Fix $2<\zeta_-<\zeta_+<\infty$. Then for any $x,y \in\mathcal{H}$ with $\zeta_-<x_\parallel,y_\parallel<\zeta_+$ and $\lvert x-y\rvert\geq1$ it holds that
    \begin{equation}
        \mathcal{G}^\mathcal{H}(x,y)\geq\frac{C}{\lvert x-y\rvert^3},
    \end{equation}
    where $C>0$ depends only on $\zeta_-,\zeta_+$.
\end{lemma}
\begin{proof}
    We have $\zeta_--2<h_x,h_y<\zeta_+-2$, hence $\lvert x^*-y\rvert\leq\lvert x-y\rvert+2h_x\leq(2\zeta_+-3)\lvert x-y\rvert$ for $\lvert x-y\rvert\geq1$. Inserting this into \eqref{eq: HalfspaceGreen} gives the claim with $C=2(\zeta_--2)^2/(\pi(2\zeta_+-3)(2\zeta_+-2))$.
\end{proof}
The lower bound of \cref{thm: ScreenDecay} now follows mutatis mutandis from the proof of the lower bound of \cref{thm: ChainDecay}: for $y\in B_{1/4}(z_i)\subset Q_i$ we have $9/4<y_\parallel<15/4$ and $5/2\leq(z_j)_\parallel\leq7/2$, while $d_{ij}-1/4\leq\lvert z_j-y\rvert\leq d_{ij}+5/4$, so that for $d_{ij}\geq5/4$ \cref{lem:hit-lower} and \cref{lem: HalfspaceGreen} (with $\zeta_-=9/4$, $\zeta_+=4$) yield $V_i(z_j)\geq c\,d_{ij}^{-3}$. For $d_{ij}<5/4$ the bound $V_i(z_j)\geq c_0'$ holds as in \cref{sssec: ChainLowerProof}, and the Harnack argument given there together with \eqref{eq: LowerFlux} transfers these bounds to $\lvert\mathcal{C}_{ij}\rvert$.
\subsection{Crystals}
It remains to prove the lower bound for crystals. In this regime, a simpler Cameron-Martin based argument is sufficient.

Once again, we proceed as in \cref{lem: Connectedness}, constructing a path $\tilde{\gamma}$ from a point $z^*\in\partial B_{r_j+d_{\min}/3}(x_j)$ to a point of $\partial B_{r_i+d_{\min}/3}(x_i)$, $i\neq j$, such that $\operatorname{dist}(\tilde{\gamma},\mathcal{D})\geq d_{\min}/3$, whose length is upper bounded by $d_{ij}$ times a constant $C$ depending on the geometric parameters, but not $i,j$, and append the radial segment to $x_i$. We parametrise $\tilde{\gamma}$ by arc length, so that $\tilde{\gamma}:[0,\ell]\to\R^3$ has unit speed and $\ell\leq Cd_{ij}+r_{\max}+d_{\min}/3$, and set $\epsilon=\min\{d_{\min}/6,r_{\min}\}$. Since under $\mathbb{P}_{z^*}$ $\{\tau<\infty, B(\tau)\in \partial D_i\}\supset \{\max_{0\leq t\leq \ell}\lvert B(t)-\tilde{\gamma}(t)\rvert <\epsilon\}$ (on this event $B$ keeps a distance greater than $d_{\min}/6$ from $\mathcal{D}\setminus D_i$ and satisfies $\lvert B(\ell)-x_i\rvert<\epsilon\leq r_i$), Cameron Martin as in the proof of \cref{lem: Killing}, now with $\int_0^\ell\lvert\tilde{\gamma}'\rvert^2dt=\ell$, yields
\begin{equation}
    V_i(z^*)\geq\exp\left(-\ell/2\right)\mathbb{P}_0(\max_{0\leq t\leq \ell}\lvert B(t)\rvert<\epsilon ).
\end{equation}
The last probability decays at most exponentially in $\ell$. Indeed, for $\lvert w\rvert\leq\epsilon/2$ the event $\{\max_{0\leq t\leq1}\lvert B(t)-(1-t)w\rvert<\epsilon/2\}$ is contained in $\{\max_{0\leq t\leq1}\lvert B(t)\rvert<\epsilon,\ \lvert B(1)\rvert<\epsilon/2\}$ under $\mathbb{P}_w$, so that Cameron Martin gives
\begin{equation}
    p_0=\inf_{\lvert w\rvert\leq\epsilon/2}\mathbb{P}_w\left(\max_{0\leq t\leq1}\lvert B(t)\rvert<\epsilon,\ \lvert B(1)\rvert<\epsilon/2\right)\geq e^{-\epsilon^2/8}\,\mathbb{P}_0\left(\max_{0\leq t\leq1}\lvert B(t)\rvert<\epsilon/2\right)>0,
\end{equation}
and the Markov property at the times $1,\dots,\lceil\ell\rceil-1$ yields $\mathbb{P}_0(\max_{0\leq t\leq\ell}\lvert B(t)\rvert<\epsilon)\geq p_0^{\lceil\ell\rceil}$. Consequently, $V_i(z^*)\geq p_0\exp\left(-\ell\left(1/2+\log(1/p_0)\right)\right)\geq c\exp(-c'd_{ij})$ with $c,c'>0$ depending only on the geometric parameters.
Applying Harnack's inequality to obtain a uniform lower bound on $\partial B_{r_j+d_{\min}/3}(x_j)$, as in the proof of the lower bound of \cref{thm: ChainDecay}, together with \cref{lem: Flux} yields the claim.

\section{Conclusion}\label{sec: Conclusion}
In this paper, we have derived upper and lower bounds off-diagonal decay estimates for the capacitance matrix of resonator systems in chain, screen and crystal geometries, using the classical connection between Brownian motion and the Dirichlet problem. Our main results, \cref{thm: ChainDecay,thm: ScreenDecay,thm: CrystalDecay}, show that the rate of decay is determined by the dimensionality of the structure.

Many interesting questions remain. For one, it is natural to expect that similar estimates hold for systems of high-index dielectric resonators \cite{ammari2023mathematical}. Another promising direction would be to derive sharper truncation estimates \cite{ammari2026resolvent}, which can be used to design fast computational methods \cite{ammari2026waveguidingsystemshighcontrast}. More ambitious, but equally important, would be decay estimates in the non-subwavelength regime, where the capacitance matrix becomes frequency dependent \cite{ammari2026frequency}.

\section*{Statements and declarations}

\subsection*{Funding}
C.~Thalhammer was supported in part by the Swiss National Science Foundation grant number 200021-236472.

\subsection*{Competing interests}
The authors have no competing interests to declare.


\subsection*{Declaration of AI use}
The authors interacted with and guided ChatGPT~6 and Claude Fable 5.1, Opus 5.5 whilst developing the proofs and then wrote the paper. The authors take full responsibility for the accuracy and integrity of all content in this paper.

\printbibliography

\end{document}